\documentclass[10pt,a4paper,twoside,notitlepage]{article}

\usepackage{amsfonts}
\usepackage{amsmath}
\usepackage{amssymb}

\usepackage{graphicx}
\usepackage{psfrag}
\usepackage{color}

\newtheorem{thm}{Theorem}
\newtheorem{lem}{Lemma}[section]

\newtheorem{prop}{Proposition}[section]
\newtheorem{defn}{Definition}[section]
\newtheorem{rem}{Remark}[section]

\newcommand{\N}{\mathbb{N}}
\newcommand{\R}{\mathbb{R}}

\newcommand{\ve}{\varepsilon}
\newcommand{\n}{\noindent}
\newcommand{\vatt}{\mathbb{E}}

\newcommand{\Div}{\mathrm{div}}
\newcommand{\Tr}{\mathrm{tr}}
\newcommand{\D}{\mathrm{D}}
\newcommand{\Id}{\mathrm{I}}
\newcommand{\Img}{\mathrm{Im}}
\newcommand{\Mat}{\mathrm{Mat}}
\newcommand{\Sym}{\mathrm{Sym_d}}
\newcommand{\Symp}{\mathrm{Sym_d^+}}

\begin{document}

\title{{\huge Linear-Quadratic $N$-person and\\ Mean-Field Games with Ergodic Cost
}}
 
\author{
  \vspace{1cm}
    {\scshape Martino Bardi
    \thanks{Dipartimento di Matematica, Universit\`a 
    di Padova,
      Via Trieste 63,
      I-35121 Padova, Italy, {\tt\small bardi@math.unipd.it}. This author  was partially supported by the Fondazione CaRiPaRo Project ``Nonlinear Partial Differential Equations: 
models, analysis, and control-theoretic problems," the MIUR project PRIN "Viscosity, geometric, and control methods for nonlinear diffusive models", and the European Project Marie Curie ITN "SADCO - Sensitivity Analysis for Deterministic Controller Design". }       
    \,  and Fabio S. Priuli 
    \thanks{Istituto per le Applicazioni del Calcolo ``M. Picone'', 
    Consiglio Nazionale delle Ricerche, 
    Via dei Taurini 19, I-00185 Roma, Italy,
    {\tt\small f.priuli@iac.cnr.it}} 
    } 
      }
             
\date{\today}

\maketitle
\begin{abstract}
We consider stochastic differential games with $N$ players, linear-Gaussian dynamics in arbitrary state-space dimension, and long-time-average cost with quadratic running cost. Admissible controls are feedbacks for which the system is ergodic.
We first study the existence of affine Nash equilibria by means of an associated system of $N$ Hamilton-Jacobi-Bellman and $N$ Kolmogorov-Fokker-Planck partial differential equations. We give necessary and sufficient conditions for the existence and uniqueness of quadratic-Gaussian solutions in terms of the solvability of suitable algebraic Riccati and Sylvester equations. Under a symmetry condition on the running costs and for nearly identical players we study the large population limit, $N$ tending to infinity, and find a unique quadratic-Gaussian solution of the pair of Mean Field Game HJB-KFP equations. 
Examples of explicit solutions are given, in particular for consensus problems. 

\end{abstract}

\section[Introduction]{Introduction} \label{sec:int}

We consider a system of $N$ stochastic differential equations
\begin{equation}\label{eq:sde}
dX^i_t=(A^iX^i_t-\alpha^i_t)dt+\sigma^idW^i_t\,, \qquad X_0^i=x^i\in\R^d\,,\qquad i=1,\ldots,N\,,
\end{equation}
where $A^i,\sigma^i$ are given $d\times d$ matrices, with $\det(\sigma^i)\neq 0$, $(W^1_t,\ldots,W^N_t)$ are $N$ independent $d$--dimensional standard Brownian motions, and $\alpha^i_t\colon [0,+\infty[\to\R^d$ is a 
process adapted to $W^i_t$ which represents the control of the $i$--th player of the differential game that we now describe. For each initial positions $X=(x^1,\ldots,x^N)\in\R^{Nd}$ we consider for the $i$--th player controls whose associated process is ergodic and the long--time--average cost functional with quadratic running cost
\begin{equation}\label{eq:cost}
J^i(X,\alpha^1,\ldots,\alpha^N):= \liminf_{T\to\infty}\,{1\over T}\,\vatt\left[
\int_0^T \,{(\alpha_t^i)^TR^i\alpha_t^i\over 2}\,+(X_t-\overline{X_i})^T Q^i(X_t-\overline{X_i})\,dt
\right]\,,
\end{equation}
where $\vatt$ denotes the expected value, $R^i$ are positive definite symmetric $d\times d$ matrices, $Q^i$ are symmetric $Nd\times Nd$ matrices, and $\overline{X_i}\in\R^{Nd}$ are given reference positions. 

 For this $N$--person game we 
study the following two problems:

\n 1. the synthesis 
of Nash equilibrium strategies in feedback form, and of the 
probability distribution for the position of each player at the equilibrium,
from a system of elliptic partial differential equations associated to the game, 

\n 2. the large population limits as $N\to\infty$ 
 of these strategies and distributions 
 and their connection with the Mean-Field Games partial differential equations introduced by Lasry and Lions~\cite{LL1, 
  LL3}. 

The first problem 
is classicaly formulated within the 
 theory 
 of Hamilton-Jacobi-Bellman equations associated to $N$-person differential games, as was done in \cite{BF1, BF} for compact state space. This leads to a system of $N$ PDEs in $\R^{Nd}$ strongly coupled in the gradients of the unknown value functions. Instead, we exploit the independence of the dynamics of different players, that makes the game merely \emph{cost-coupled}, and follow the approach of Lasry and Lions~\cite{LL1, 
  LL3} leading to a system of $N$ nonlinear PDEs in $\R^{d}$ of HJB  type coupled with $N$ Kolmogorov-Fokker-Planck equations for the invariant measure of the process associated to the Nash equilibrium. 
  This has several advantages, 
  including a much weaker coupling of the new system of PDEs. 
We refer to \cite{BardiFeleqi} for more information on the connections between the two approaches. In view of the Linear-Quadratic-Gaussian structure of the game we look for solutions of the HJB-KFP system in the class of quadratic value functions and multivariate Gaussian distributions and give necessary and sufficient conditions for both existence and uniqueness of solutions in this class. This produces Nash equilibria in the form of affine feedbacks.

The second problem is set in the framework of \emph{nearly identical players}, as in \cite{Bardi}. We first characterize the existence and uniqueness of identically distributed solutions, then take the limit of these solutions as $N\to\infty$ and show that it solves the system of two Mean-Field Games PDEs in $\R^d$
\begin{equation}\label{eq:hjkfp_limit*}
\left\{
\begin{array}{l}
-\Tr\left(\frac{\sigma\sigma^T}2\D^2 v\right)+\displaystyle\,{1\over 2}\,\nabla v^T R^{-1}\nabla v-\nabla v^TAx+\lambda=\hat V[m](x)\\
-\Tr\left(\frac{\sigma\sigma^T}2\D^2 m\right)-\displaystyle\Div\left(m \cdot(R^{-1}\nabla v-A x)\right)=0\\
\int_{\R^d}m(x)\,dx=1\,,\qquad m>0\,,
\end{array}
\right.
\end{equation}
in the unknowns $(v,m,\lambda)$, with $v, m\in C^2(\R^d)$
 and $\lambda\in \R$, where $\Tr$ and $\Div$ are 
the trace of a matrix and the divergence operator, respectively, and $\hat V[m]$ is an integral operator sending probability densities into quadratic polynomials defined in terms of the blocks of the matrix $Q^i$. 
Moreover, the solution obtained in the large population limit is the unique having $v$ quadratic and $m$ Gaussian.
We also show that such solution is unique among 
general solutions of \eqref{eq:hjkfp_limit*} under a condition on a submatrix of $Q^i$ meaning that \emph{imitation is not rewarding} in the large population limit. 

The strategy 
of proof for  both problems is 
 the same: we insert the coefficients of the quadratic value function and of the Gaussian distribution 
 into the system of partial differential equations and we show that this reduces the problem to the solvability of an algebraic Riccati equation and a Sylvester equation. Besides proving the existence of solutions these matrix equations can be used to solve explicitly some examples or can be solved numerically in more complex cases. There is a large literature on the numerical resolution of Riccati equations for which we refer, e.g., to \cite{BIM} and 
 its bibliography. 
 
We find explicit formulas for the solutions of \eqref{eq:hjkfp_limit*} in the case the diffusion matrix $\sigma$ and the cost matrix $R$ are constants times the identity matrix, and the drift $A$ of the system is either a symmetric or a non-defective matrix. 
In particular we treat 
 cost functionals depending on the state via the quadratic form
\begin{equation}\label{eq:special_cost*}
F^i(X^1,\ldots,X^N)=\frac 1{N-1}\sum_{j\neq i} (X^i-X^j)^T P^N(X^i-X^j) , \qquad i=1,\dots, N ,
\end{equation}
with $P^N\to \hat P>0$, a model arising in \emph{consensus problems}.
We recall that a consensus process aims at reaching
an agreement among several agents 
on some common state properties. This is an active area of research within multi-agent control and coordination.
We refer to \cite{NCMH} for the motivations of consensus problems, their connection with mean field control theory, and references to
the large literature on the subject. Note that the cost functional \eqref{eq:special_cost*} pushes the players to take positions close to each other. We show the existence of a quadratic-Gaussian solution with mean $\mu$ for each solution of $A^i\mu=0$. Therefore such solutions and the corresponding Nash equilibria can be infinitely many: note that here imitation is rewarding.

Most of the results of this paper and explicit formulas for the solutions were derived by the first-named author in \cite{Bardi} in the case of $1$--dimensional state space, i.e., $d=1$, where the analysis is much simpler because the search for quadratic-Gaussian solutions leads to scalar polynomial equations of degree 
at most two. For $d>1$, instead, we arrive at some nontrivial 
algebraic Riccati equations coupled with Sylvester equations that require a much heavier use of matrix algebra and do not admit explicit solutions in general. 
In a sequel of the present paper~\cite{Priuli} the second-named author studies several singular limits of the $N$-person and Mean-Field games considered here, such as the vanishing viscosity, the cheap control, and the vanishing discount limit.



Linear-Quadratic differential games have a large literature, see the books \cite{BO, Engw} and the references therein.
Large population limits for multi-agent systems 
 were studied by Huang, Caines and Malhame, independently of Lasry-Lions. They introduced a method named Nash certainty equivalence principle \cite{HCM:03, HCM:06, HCM:07ieee} that first produces a feedback from a  mean-field equation and then shows that it is an $\epsilon$-Nash equilibrium for the $N$-person game if $N$ is large enough.
We cannot review here the number of papers inspired by their approach, but let us mention 
 \cite{LiZ} and 
  \cite{BSYY} 
   for LQ problems, \cite{BauTBas} about robust control, 
    \cite{NC} for recent progress on nonlinear systems, \cite {NCMH} for consensus problems, \cite{KTY} on the rate of convergence as $N\to\infty$, and the references therein. Some of these papers also deal with ergodic cost functionals, e.g., \cite{HCM:07ieee, LiZ}, but their assumptions and methods differ from ours.

Concerning the Lasry-Lions approach to MFG, besides their pioneering papers~\cite{LL1, LL2, LL3}
 let us mention the lecture notes \cite{Car} and \cite{GLL}, \cite{ACD} on numerical methods, \cite{GMS} on discrete games, \cite{CLLP} on the long time behaviour of solutions, \cite{Fel} on the large population limit for nonlinear ergodic control of several populations, and the thesis \cite{Cir} on multi-population models.
A very recent survey 
 on MFG focusing on the comparison with the 
  theory of 
 mean-field type control 
 is \cite{BFY}.
There  is a wide spectrum of applications of Mean-Field Games that we do not try to list here 
 and refer instead to the quoted literature.

The paper is organised as follows.
\n In Section~\ref{sec:prelim} we define admissible strategies, introduce the system of HJB--KFP PDEs associated to the $N$--person game and recall some known facts about matrices and algebraic Riccati equations.
\n In Section~\ref{sec:n_play} we present our main result about Nash equilibrium strategies for the $N$--players game. 
\n In Section~\ref{sec:nearly_id} we define the games with nearly identical players and give the existence and uniqueness result in that case.
\n Section~\ref{sec:limit} is devoted to the analysis of the limit when the number of players tends to infinity, under natural rescaling assumptions on the matrix coefficients of the game. 
\n Finally, in Section~\ref{sec:explicit} we present various explicit sufficient conditions for the validity of the previous theory and some explicit solutions, in particular for consensus problems.


Some selected results of this paper were presented at the 52nd IEEE-CDC in Florence, December 2013 \cite{BardiPriuli:cdc}.

\section[Preliminaries]{Preliminaries} \label{sec:prelim}

\subsection{Some properties of symmetric matrices}

In the following, we will use the notation $\Mat_{d\times d}(\R)$ for the linear space of real $d\times d$ matrices, $\Sym$ for the subspace of real symmetric $d\times d$ matrices and $\Id_d$ for the identical $d\times d$ matrix.

Given a matrix $M\in\Sym$, we say that $M$ is positive semidefinite (resp. positive definite) if for all $x\in\R^d$ there holds $x^TMx\geq 0$ (resp. if for all $x\in\R^d\setminus\{0\}$ there holds $x^TMx> 0$). We will also use the notation $\Symp$ for the set of real symmetric and positive definite $d\times d$ matrices. 

Finally, given a (real) $d\times d$ matrix $M$, we denote its spectrum by $\mathrm{spec}(M)$. Recall that eigenvalues of a matrix $M$ depend continuously on the coefficients of the matrix (see~\cite{SerreMat}), so that if we have $M_n\to M$, then in particular the eigenvalues of $M_n$ converge to eigenvalues of $M$.

We summarize in the following proposition some facts that will be used to prove our results. 

\begin{prop}\label{prop:matrices} The following facts hold.
\begin{description}
\item{\it (i)} Let $H\in\Sym$ and $K\in\Symp$. Then, $HK$ is diagonalizable with real eigenvalues. Moreover, the number of positive (resp. negative) eigenvalues of $HK$ is equal to the number of positive (resp. negative) eigenvalues of $H$. The same holds for $KH$.
\item{\it (ii)} By defining for every matrix $M\in\Sym$
\begin{equation}\label{eq:max_spec_norm}
\|M\|:=\max\,\{|\ell|~;~\ell\in\mathrm{spec}(M)\}\,,
\end{equation}
we obtain a norm. Moreover, if $M$ is positive semidefinite, then $\|M\|$ is simply the largest eigenvalue of $M$, i.e.
$$
\|M\|=\max\mathrm{spec}(M)\,.
$$
\item{\it (iii)} Let $H,K\in\Sym$ be positive semidefinite matrices and let $L\in\Symp$. If $H-K$ is positive semidefinite, then $\|H\|\geq\|K\|$ and 
$$
\max\mathrm{spec}(HL)\geq \max\mathrm{spec}(KL)\,.
$$
\end{description}
\end{prop}

\n{\it Sketch of the proof.} We refer to~\cite{SerreMat} for a proof of {\it (i)} and {\it (ii)}. Concerning {\it (iii)}, the part about $\|H\|\geq\|K\|$ follows from Weyl's inequalities (see~\cite{SerreMat}): by setting $\lambda^*(M):=\max\mathrm{spec}(M)$, we obtain
$$
\|K\|=\lambda^*(K)\leq\lambda^*(H)+\lambda^*(K-H)\leq \lambda^*(H)=\|H\|\,.
$$
The remaining part requires to recall that the spectrum of a matrix remains the same under changes of basis, and that $H-K\geq 0$ implies $M^THM - M^TKM\geq 0$ for all matrices $M\in\Mat_{d\times d}(\R)$.
Therefore, we obtain
$$
\lambda^*(HL)=\lambda^*(\sqrt{L}HL(\sqrt{L})^{-1})=\lambda^*(\sqrt{L}H\sqrt{L})
\geq\lambda^*(\sqrt{L}K\sqrt{L})=
\lambda^*(KL)\,,
$$
and this completes the proof.~~$\diamond$

\subsection{Admissible strategies}

We consider strategies whose corresponding solution to~\eqref{eq:sde} is ergodic.

\begin{defn}\label{defn:admiss_strategy} A strategy $\alpha^i_t$ is said to be \emph{admissible} (for the $i$--th player) if it is a 
 process adapted to the Brownian motion $W^i_t$, with $\vatt[|\alpha^i_t|^2]$ bounded on $[0,T]$ for all $T$, and  such that the corresponding solution $X^i_t$ to~\eqref{eq:sde} satisfies
\begin{itemize}
\item 
$\vatt[
(X^i_t)(X^i_t)^T]$
 is bounded on $[0,T]$ for every $T$;
\item $X^i_t$ is \emph{ergodic} in the following sense: there exists a probability measure $m^i=m^i(\alpha^i)$ on $\R^d$ such that
$$
\int_{\R^d}|x|\,dm^i(x)<\infty
\qquad\qquad
\int_{\R^d}|x|^2\,dm^i(x)<\infty
$$
and
$$
\lim_{T\to+\infty}\,{1\over T}~\vatt\left[\int_{\R^d}g(X^i_t)\,dt\right]=\int_{\R^d}g(x)\,dm^i(x)\,,
$$
locally uniformly w.r.t. the initial state $X^i_0$, for all functions $g$ which are polynomials of degree at most $2$.
\end{itemize}
\end{defn}

\n Below, we prove that affine strategies are admissible. Preliminarily, let us introduce the notation ${\cal N}(\mu,V)$ for a multivariate Gaussian distribution with mean $\mu\in\R^d$ and covariance $V\in\Symp$, i.e. a distribution with density
$$
\gamma\exp\left\{-\,{1\over 2} (x-\mu)^TV^{-1}(x-\mu)\right\}\,,\qquad\mbox{where}\qquad \gamma:=(2\pi)^{-d/2}\sqrt{\mathrm{det}(V)}\,.
$$
Moreover, we say that a probability measure is an \emph{invariant measure} for the process $X_t$ in $\R^d$ if 
\[
\int_{\R^d}\vatt[\phi(X_t)\, |\, X_0=y]\, d \mu(y) = \int_{\R^d} \phi(x)\, d \mu(x)
\]
for all $t\geq 0$ and $\phi : \R^d\to\R$ bounded and uniformly continuous.

\begin{prop}\label{prop:affine_admissible} For the affine feedback
\begin{equation}\label{eq:aff_feedback}
\alpha^i(x)= K^ix+c^i\,,\qquad\qquad x\in\R^d\,,
\end{equation}
with $K^i\in\Mat_{d\times d}(\R)$ such that the matrix $A^i-K^i$ has only eigenvalues with negative real part, and $c^i\in\R^d$, consider the process $\alpha^i_t:= \alpha^i(X^i_t)$ where $X^i_t$ solves
\begin{equation}\label{eq:csde}
dX^i_t=[(A^i-K^i)X^i_t-c^i]dt+\sigma^idW^i_t\,.
\end{equation}
Then $\alpha^i_t$ is admissible. Moreover, the process $X^i_t$ has a unique invariant measure $m^i$ given by a multivariate Gaussian ${\cal N}(\mu,V)$ with mean $\mu=(A^i-K^i)^{-1}c^i$ and covariance matrix $V$ which satisfies the algebraic relation
\begin{equation}\label{eq:var_equation}
(A^i-K^i) V+V(A^i-K^i)^T+\sigma^i(\sigma^i)^T=0\,,
\end{equation}
and $X^i_t$ is ergodic w.r.t. such a measure $m^i$.
\end{prop}

\begin{rem}\label{rem:var_equation} \rm{
Note that the equation~\eqref{eq:var_equation} satisfied by the covariance matrix $V$ admits a unique solution. Indeed, it is a Sylvester equation of the form
$$
M Y + Y N= -\sigma^i(\sigma^i)^T\,,
$$
with $M=N^T$, so that the matrices $M$ and $-N$ have no eigenvalues in common (see~\cite[Lemma~2.31]{Engw}).
}\end{rem}


\n {\it Proof of Proposition~\ref{prop:affine_admissible}.} It is well known (see for instance~\cite[Theorem~4.3]{BressanPiccoli}) that given a matrix $M\in\Mat_{d\times d}(\R)$ whose eigenvalues have all strictly negative real part, it is always possible to find a matrix $P\in\Symp$ such that $M^T P+PM=-\Id$. 
\n We now claim that, by denoting with $P^i$ the symmetric positive definite matrix corresponding to $M=A^i-K^i$, $V(x):=
 x^T P^i x$ is a Lyapunov--like function for the system~\eqref{eq:csde}. Indeed, by denoting with ${\cal L}$ the infinitesimal generator of the process, we have
$$
{\cal L}V(x)=\Tr\left({\sigma^i(\sigma^i)^T\over 2} P^i\right) + x^TP^i (A^i-K^i)\,x + \big((A^i-K^i)\,x\big)^T P^ix=\Tr\left({\sigma^i(\sigma^i)^T\over 2} \,P^i\right) - |x|^2\,,
$$
which is strictly negative outside a ball of radius $R>\sqrt{{\sigma^i(\sigma^i)^T\over 2} \,P^i}$. Hence, the existence of a unique invariant measure $m^i$ for~\eqref{eq:csde} follows by exploiting the theory of Khasminskii~\cite{Khas} or the results in~\cite{BardiCesaroniManca}. Observing that if a unique invariant measure exists, then such a measure is also ergodic in the sense of Definition~\ref{defn:admiss_strategy} (see e.g.~\cite[Theorem~5.16]{DaPrato}), we have verified the second property required by admissibility.

\smallskip

It is well known (cf.~\cite[Section~5.6]{KarShreve}) that, 
for solutions of linear stochastic equations~\eqref{eq:csde}, 
the mean vector ${\frak m}_i(t):=\vatt[X^i_t]$ and the covariance matrix ${\frak v}_i(t):=\vatt[(X^i_t)(X^i_t)^T]$ are respectively solutions of
$$
\dot {\frak m}_i(t)= (A^i-K^i)\, {\frak m}_i(t) -c^i\,,
$$
and
$$
\dot {\frak v}_i(t)=(A^i-K^i)\, {\frak v}_i(t)+{\frak v}_i(t)(A^i-K^i)^T+\sigma^i(\sigma^i)^T\,,
$$
whence boundedness of first and second moments follows.

Finally, since a multivariate Gaussian ${\cal N}(\mu, V)$ with $\mu=(A^i-K^i)^{-1}c^i$ and $V$ solving~\eqref{eq:var_equation} is indeed a stationary solution of~\eqref{eq:csde},
by uniqueness we get $m^i={\cal N}(\mu, V)$.~~$\diamond$


\subsection{Algebraic Riccati equations}

We recall here some basic facts about algebraic Riccati equations (ARE in the following).

\begin{prop}\label{prop:ARE} Consider the ARE
\begin{equation}\label{eq:ARE}
Y {\cal R} Y
- {\cal Q}=0
\end{equation}
with ${\cal R}\in\Symp$ and ${\cal Q}\in\Sym$, and introduce the following notations
\begin{equation}\label{eq:inv_spaces}
\Xi_S:= \left[
\begin{array}{c}
I_d\\ S
\end{array}
\right]\in\Mat_{2d\times d}(\R)\,,
\qquad\qquad
{\cal H}:=
\left(
\begin{array}{cc}
{\bf 0} & {\cal R}\\
{\cal Q} & {\bf 0}
\end{array}
\right)\in\Mat_{2d\times 2d}(\R)\,,
\end{equation}
where $S$ is any element of $\Mat_{d\times d}(\R)$, and $\Img\,\Xi_S$ for the $d$--dimensional linear subspace of $\R^{2d}$ spanned by the columns of $\Xi_S$. Then the following facts hold.
\begin{description}
\item{\it (i)} $Y$ is a solution of~\eqref{eq:ARE} if and only if $\Img\,\Xi_Y$ is ${\cal H}$--invariant, i.e. if and only if ${\cal H}\xi\in\Img\,\Xi$ for all $\xi\in\Img\,\Xi$.
\item{\it (ii)} If the matrix ${\cal H}$ has no purely imaginary nonzero eigenvalues, then equation~\eqref{eq:ARE} has solutions $Y$ such that $Y=Y^T$.
\item{\it (iii)} If ${\cal Q}$ is positive definite, then all eigenvalues of ${\cal H}$ are real and different from zero, and equation~\eqref{eq:ARE} has a unique symmetric solution $Y$ such that the matrix ${\cal R}Y$ has only positive eigenvalues. In particular, 
\begin{equation}\label{eq:spec_ARE}
\mathrm{spec}({\cal R}Y) = \mathrm{spec}({\cal H})\cap (0,+\infty)\,,
\end{equation}
and $Y$ is also the unique symmetric positive definite solution to~\eqref{eq:ARE}.
\end{description}
\end{prop}

\n{\it Sketch of the proof.} The proof follows from standard arguments about Riccati equations that can be found in~\cite{Engw,LR}. We give here some explicit references for sake of completeness. Part {\it (i)} is contained in Proposition 7.1.1 of~\cite{LR}. Part {\it (ii)} is a particular case of Theorem 8.1.7 in~\cite{LR}. Concerning {\it (iii)}, if we assume ${\cal Q}$ positive definite, then $\ell$ is an eigenvalue for ${\cal H}$ if and only if $\ell$ is a solution of the equation
$$
0=\det({\cal H}-\ell I_{2d})=
\det\left(
\begin{array}{cc}
-\ell I_d & {\cal R}\\
{\cal Q} & -\ell I_d
\end{array}
\right)=\det(\ell^2 I_d- {\cal R}{\cal Q})
$$
i.e. if and only if $\ell^2$ is an eigenvalue of the $d\times d$ matrix ${\cal R}{\cal Q}$. But we are assuming that  both ${\cal R}$ and ${\cal Q}$ are positive definite, so ${\cal R}{\cal Q}$ has only positive eigenvalues and therefore all eigenvalues of ${\cal H}$ are in $\R\setminus\{0\}$. Now Theorem 8.3.2 of~\cite{LR} ensures that there is a unique solution $Y$ of~\eqref{eq:ARE} in $\Sym$ such that ${\cal R}Y$ has only positive eigenvalues and such solution is characterized by~\eqref{eq:spec_ARE}. In turn,~\eqref{eq:spec_ARE} implies that $Y\in\Symp$. Finally, assume there is another solution $Z\in\Symp$, then ${\cal R}Z$ would have only positive eigenvalues too, and hence we would have $Z=Y$ by the characterization above.~~$\diamond$

\subsection{HJB and KFP equations associated to the $N$--person game}

\n We focus our attention on the game~\eqref{eq:sde}--\eqref{eq:cost}, in order to write the system of HJB--KFP equations associated to the game~\eqref{eq:sde}--\eqref{eq:cost}, as in~\cite{Bardi, LL1, LL3}.
We start by remarking that the part of the cost depending on the state of the game can be also written as
\begin{align}\label{eq:quad_cost}
F^i(X^1,\ldots, X^N)&:= (X-\overline{X_i})^T Q^i(X-\overline{X_i})
=\sum_{j,k=1}^N (X^j-\overline{X_i}^j)^T Q^i_{jk}(X^k-\overline{X_i}^k)
\end{align}
where the matrices $Q^i_{jk}$ are $d\times d$ blocks of $Q^i$. The standing assumptions on the game are summed up in the following condition.

\begin{description}
\item{\bf (H)} Assume that $\sigma^i$ in~\eqref{eq:sde} are invertible matrices, that $R^i$ in~\eqref{eq:cost} belong to $\Symp$ and that $Q^i$ in~\eqref{eq:cost} are symmetric matrices. Moreover, assume that blocks $Q^i_{ii}$ belong to $\Symp$ as well. 
\end{description}

\begin{rem} \rm{
Each block $Q^i_{jk}$ collects the costs for player $i$ per unit of displacement of players $j$ and $k$ from their reference positions $\overline{X_i}^j$ and $\overline{X_i}^k$, respectively. The assumption 
 $Q^i_{ii}>0$ means that his own reference position $\overline{X_i}^i$ is a preferred position for the $i$--th player. This condition has a clear interpretation and makes some calculations easier, but for the validity of Theorem~\ref{thm:Nplay} it can be weakened to $Q^i_{ii}+(A^i)^TR^iA^i/2>0$.
}\end{rem}

For the game~\eqref{eq:sde}--\eqref{eq:cost} under consideration, we observe that the $i$-th Hamiltonian takes the form
\begin{align*}
H^i(x,p)&:=\min_{\omega}\left\{-\omega^T \,{R^i\over 2}\,\omega - p^T\big(A^ix-\omega\big)\right\}=-p^TA^ix+\min_{\omega}\left\{-\omega^T \,{R^i\over 2}\,\omega - p^T
\omega\right\}\,.
\end{align*}
Since the minimum is attained at $(R^i)^{-1} p$, we get 
\begin{align*}
H^i(x,p)&=-((R^i)^{-1} p)^T \,{R^i\over 2}\,((R^i)^{-1} p) - p^T\big(A^ix-(R^i)^{-1} p\big)=p^T\,{(R^i)^{-1}\over 2}\,p - p^TA^ix\,.
\end{align*}
we introduce  the notations
\begin{equation}\label{eq:rhs}
f^i(x;{\frak m}^1,\ldots,{\frak m}^N):=
\int_{\R^{d(N-1)}} F^i(\xi^1,\ldots,\xi^{i-1},x,\xi^{i+1},\ldots\xi^N)\prod_{j\neq i} d{\frak m}^j(\xi^j)\,,
\end{equation}
for any $N$--vector of probability measures $({\frak m}^1,\ldots,{\frak m}^N)$, and
$$
\nu^i:=\,{ (\sigma^i)(\sigma^i)^T\over 2}\in\Mat_{d\times d}(\R)\,.
$$
The classical system of $N$ Hamilton-Jacobi-Bellman equations in $\R^{Nd}$ whose solutions generate Nash feedback equilibria  is,
for $i=1,\ldots,N$,
\begin{equation}\label{eq:hjk}
-\sum_{j=1}^N\Tr(\nu^i\,\D^2_{x^j} v^i)+H^i(x^i,\nabla_{x^i} v^i)+\sum_{j\ne i}\frac{\partial H^j}{\partial p}(x^j,\nabla_{x^j} v^j)\cdot \nabla_{x^j} v^i+\lambda^i=F^i(X) , 
\end{equation}
see \cite{BF1, BF} and the references therein, where $X:=(x^1,\dots,x^N)$, and  the unknowns are $v^i(X)$ and the constants $\lambda^i$, $i=1,\dots,N$. Note that this system of PDEs is strongly coupled via the terms $\nabla_{x^j} v^j$ appearing in the $i$-th equation.
In view of the independence of the dynamics of the different players 
 we follow the approach by Lasry and Lions \cite{LL1, LL3} and suppose the solution $v^i$ of \eqref{eq:hjk} depends only on $x^i,$ $i=1,\dots,N$. We consider the linearization  around $v^i$ of the $i$-th equation of \eqref{eq:hjk}, that is,
 \[
 -\Tr(\nu^i\,\D^2
  m^i)-\displaystyle\Div\left(m^i {\partial H^i\over\partial p}(x
 ,\nabla v^i)\right)=0 , \quad \text{in } \R^d , \;i=1,\dots,N ,
\]
and assume this system has positive solutions $m^i$ with $\int_{\R^d}m^i(x)\,dx=1$. Next we multiply the $i$-th equation in \eqref{eq:hjk} by $\Pi_{j\ne i} m^j(x^j)$ and integrate over $\R^{(N-1)d}$ with respect to $dx^j$, ${j\ne i}$.
Then we arrive at the system of HJB--KFP equations in $\R^d$ 
\begin{equation}\label{eq:hjkfp}
\left\{
\begin{array}{l}
-\Tr(\nu^i\,\D^2 v^i)+H^i(x,\nabla v^i)+\lambda^i=f^i(x;m^1,\ldots,m^N)\\
~\\
-\Tr(\nu^i\,\D^2 m^i)-\displaystyle\Div\Big(m^i {\partial H^i\over\partial p}(x,\nabla v^i)\Big)=0\\
~\\
\int_{\R^d}m^i(x)\,dx=1\,,\qquad m^i>0 ,
\end{array}
\right.
\qquad\qquad i=1,\ldots, N
\end{equation}
where $v^i:\R^d\to\R$
, $\lambda^i$ are real numbers and, with a slight abuse of notations, we have denoted with $m^i$ a measure as well as its density. As in~\cite{Bardi}, since we are not in the periodic setting of~\cite{LL1,LL3}, the solutions $v^i$ are expected to be unbounded and cannot be normalized by prescribing the value of their average. 
This is the system we will study in the next two sections. It will be fully justified by Theorem \ref{thm:Nplay}, where we construct directly a Nash equilibrium in feedback form from a solution of \eqref{eq:hjkfp} without resorting to the harder system \eqref{eq:hjk}.

\section[$N$--person game]{$N$--person game} \label{sec:n_play}

First of all, define the following auxiliary matrix ${\cal B}\in\Mat_{Nd\times Nd}(\R)$ as 
\begin{equation}\label{eq:matrixB}
{\cal B}:=\big({\cal B}_{\alpha\beta}\big)_{\alpha,\beta=1,\ldots,N}\,,\qquad
{\cal B}_{\alpha\beta}:= -Q^\alpha_{\alpha\beta}-\delta_{\alpha\beta}\,{(A^\alpha)^T R^\alpha A^\alpha\over 2}\,\in\Mat_{d\times d}(\R)\,,
\end{equation}
where $\delta_{\alpha\beta}$ is the Kronecker delta, and an auxiliary vector 
\begin{equation}\label{eq:vecP}
P:= \left(
\begin{array}{c}
-\sum_{j=1}^N Q^1_{1j}\,\overline{X_1}^j\\
\vdots\\
-\sum_{j=1}^N Q^N_{Nj}\,\overline{X_N}^j
\end{array}
\right)\in\R^{Nd}\,.
\end{equation}
Also, denote with $[{\cal B}, P]\in\Mat_{Nd\times (Nd+1)}(\R)$ the matrix whose columns are the columns of ${\cal B}$ and the vector $P$, i.e. 
\begin{equation}\label{eq:expanded_matrix}
[{\cal B}, P]:=\left({\cal B}^1,\ldots,{\cal B}^{Nd}, P\right)\,,
\end{equation}
being ${\cal B}^j$ the columns of the matrix ${\cal B}$.
With these notations we can state the following 
 conditions for existence and uniqueness of solution to the associated system of Hamilton--Jacobi--Bellman and Kolmogorov--Fokker--Planck equations. 
\begin{description}
\item{\bf (E)} For each $i\in\{1,\ldots,N\}$, every symmetric and positive definite solution $Y$ of the algebraic Riccati equation
\begin{equation}\label{eq:riccati}
Y\,{\nu^i R^i\nu^i\over 2}\,Y=\,{(A^i)^T R^iA^i\over 2}\,+Q^i_{ii}\,,
\end{equation}
is also a solution to the Sylvester equation
\begin{equation}\label{eq:sylvester}
Y \nu^iR^i-R^i\nu^i Y= R^i A^i-(A^i)^TR^i\,.
\end{equation}
Moreover, the matrices ${\cal B}\in\Mat_{Nd\times Nd}(\R)$ and $[{\cal B}, P]\in\Mat_{Nd\times (Nd+1)}(\R)$ have the same rank, where ${\cal B}$ is the matrix defined in~\eqref{eq:matrixB}, $P$ is the vector defined in~\eqref{eq:vecP}, and $[{\cal B}, P]$ is the matrix defined in~\eqref{eq:expanded_matrix}.
\end{description}

\begin{description}
\item{\bf (U)} The block matrix ${\cal B}$ defined by~\eqref{eq:matrixB} is invertible.
\end{description}


Some explicit conditions on the data ensuring {\bf (E)} are discussed in Section \ref{sec:explicit}. Our main result for the games with $N$ players is the following.

\begin{thm}\label{thm:Nplay} Assume that the $N$--players game having dynamics~\eqref{eq:sde} and costs~\eqref{eq:cost} satisfies assumptions {\bf (H)}. Then, the associated system of $2N$ HJB--KFP equations~\eqref{eq:hjkfp} admits solutions $(v^i,m^i,\lambda^i)_{1\leq i\leq N}$ of the form $v^i$ quadratic function and $m^i$ multivariate Gaussian, i.e. of the form
\begin{equation}\label{eq:quad_gauss}
v^i(x)=x^T\,{\Lambda^i\over 2}\, x+(\rho^i)^T x\,,\qquad\qquad
m^i={\cal N}(\mu^i,(\Sigma^i)^{-1})\,,\qquad\qquad
\lambda^i\in\R\,,
\end{equation}
for suitable symmetric matrices $\Lambda^i$ and $\Sigma^i$, $\Sigma^i$ positive definite, and vectors $\rho^i,\mu^i\in\R^d$, if and only if condition {\bf (E)} is satisfied. 
In particular, $\Sigma^i$ is the unique solution in $\Symp$ of the Riccati equation \eqref{eq:riccati}, 
$\Lambda^i=R^i(\nu^i\Sigma^i+A^i)$,  $\mu=(\mu^1,\dots,\mu^N)$ solves ${\cal B}\mu =P$, and $\rho^i=-R^i\nu^i\Sigma^i\mu^i$.

\n Moreover, solutions of the form~\eqref{eq:quad_gauss} are unique if and only if condition {\bf (U)} is satisfied and, if this is the case, the affine feedbacks
\begin{equation}\label{eq:feedback_strateg}
\overline{\alpha}^i(x)=(R^i)^{-1}\nabla v^i(x)\,,\qquad\qquad x\in\R^d,~~i=1,\ldots,N
\end{equation}
provide a Nash equilibrium strategy for all initial positions $X\in\R^{Nd}$, among the admissible strategies, and $J^i(X,\overline{\alpha})=\lambda^i$ for all $X$ and all $i$.
\end{thm}


\n {\it Proof.} The proof will be divided in several steps. Steps~1 to 4 study the particular form taken by system~\eqref{eq:hjkfp} when solutions are assumed to have the form~\eqref{eq:quad_gauss}. Step~5 proves the equivalence between existence and {\bf (E)}. Step~6 proves the equivalence between uniqueness and {\bf (U)}. Finally, Step~7 shows that the affine strategies~\eqref{eq:feedback_strateg} give a Nash equilibrium for the game~\eqref{eq:sde}--\eqref{eq:cost}.


\n {\bf Step~1.} We start by inserting functions of the form~\eqref{eq:quad_gauss} into the system~\eqref{eq:hjkfp}, 
beginning with the second equation (KFP). Notice that the hypothesis on the measure $m^i$ can be rewritten in terms of its density, which we denote again as $m^i$, as follows
\begin{equation}\label{eq:ansatz_meas}
m^i(x)=\gamma^i\exp\left\{-\,{1\over 2} (x-\mu^i)^T\Sigma^i(x-\mu^i)\right\}\,,\quad \gamma^i=(2\pi)^{-d/2}\sqrt{\mathrm{det}(\Sigma^i)}\,,
\end{equation} 
for a matrix $\Sigma^i\in\Symp$ and a vector $\mu^i\in\R^d$. 
In particular, from~\eqref{eq:ansatz_meas} we deduce $\nabla m^i(x)=-m^i(x)\Sigma^i(x-\mu^i)$.
Similarly, the condition on the value function can be rewritten in terms of its gradient as follows
\begin{equation}\label{eq:ansatz_value}
\nabla v^i(x) = \Lambda^i x+\rho^i\,,
\end{equation}
Hence, by substituting these expressions in the second equation of~\eqref{eq:hjkfp} and recalling the expressions of $H^i$ and $\nabla m^i$, we obtain
\begin{align*}
0&=-\,\Tr(\nu^i\,\D^2 m^i)-\displaystyle\Div\Big(m^i {\partial H^i\over\partial p}(x,\nabla v^i)\Big)
=-\,\Div\Big(\nu^i\nabla m^i+m^i {\partial H^i\over\partial p}(x,\nabla v^i)\Big)\\
&=\Div\Big(m^i\big( \nu^i\Sigma^i(x-\mu^i)-(R^i)^{-1}\nabla v^i+A^i x \big) \Big)\,.
\end{align*}
Therefore,~\eqref{eq:ansatz_value} implies
\begin{align*}
0&=\Div\Big(m^i\big( \nu^i\Sigma^i(x-\mu^i)-(R^i)^{-1}\Lambda^i x-(R^i)^{-1}\rho^i+A^i x \big) \Big)\\
&=m^i\,\Tr(\nu^i\Sigma^i-(R^i)^{-1}\Lambda^i+A^i)\\
&~~~~~~~~~~~~~~~~~~~~~+\nabla m^i\cdot\big[ \nu^i\Sigma^i(x-\mu^i) +(A^i-(R^i)^{-1}\Lambda^i)\,x-(R^i)^{-1}\rho^i\big]\\
&=-\,m^i~\Big\{ (x-\mu^i)^T\Sigma^i\nu^i\Sigma^i(x-\mu^i)+ (x-\mu^i)^T\Sigma^i(A^i-(R^i)^{-1}\Lambda^i)\,x\\
&~~~~~~~~~~~~~~~~~~~~~-(x-\mu^i)^T\Sigma^i(R^i)^{-1}\rho^i -\,\Tr(\nu^i\Sigma^i-(R^i)^{-1}\Lambda^i+A^i)\Big\}
\end{align*}
Since $m^i(x)>0$, this means that the other factor must vanish for every $x\in\R^d$. But such a factor is a quadratic form, and this means that its coefficients must be zero. In turn, this leads to the following matrix relations
\begin{equation}\label{eq:KFP_helper}
\nu^i\Sigma^i-(R^i)^{-1}\Lambda^i+A^i=0\,,\qquad\qquad \nu^i\Sigma^i\mu^i+(R^i)^{-1}\rho^i=0\,.
\end{equation}

\n In conclusion, necessary and sufficient condition for having solutions to the KFP equation, of the form~\eqref{eq:quad_gauss}, is that the value function $v^i$ is related to the measure $m^i$ through
\begin{equation}\label{eq:KFP_matrix}
\Lambda^i=R^i\big(\nu^i\Sigma^i+A^i\big)\,,\qquad\qquad \rho^i=-R^i\nu^i\Sigma^i\mu^i\,,
\end{equation}
where the equality for $\Lambda^i$ also imposes that the matrix $R^i\big(\nu^i\Sigma^i+A^i\big)$ is symmetric.


\n {\bf Step~2.} Let us consider now the first equation of~\eqref{eq:hjkfp} (HJB). By exploiting~\eqref{eq:KFP_matrix}, we have
$$
\Tr(\nu^i\,\D^2 v^i)=\Div(\nu^i\,\nabla v^i)=
\Div(\nu^iR^i\nu^i \Sigma^i (x-\mu^i)+\nu^iR^iA^i x)=
\Tr(\nu^iR^i\nu^i\Sigma^i+\nu^iR^iA^i)\,,
$$
$$
H^i(x,\nabla v^i)=\,{1\over 2}(\nabla v^i)^T(R^i)^{-1}\nabla v^i-(\nabla v^i)^TA^i x\,,
$$
hence the equation can be rewritten as
\begin{align*}
-\Tr(\nu^iR^i\nu^i\Sigma^i+\nu^iR^iA^i)&+\,{1\over 2}\, (x-\mu^i)^T\Sigma^i\nu^i R^i\nu^i\Sigma^i(x-\mu^i)
+\,{1\over 2}\, x^T(A^i)^T R^iA^ix\\
&+\,{1\over 2}\Big(x^T  (A^i)^T R^i\nu^i\Sigma^i (x-\mu^i)+(x-\mu^i)^T \Sigma^i\nu^i R^i A^i x\Big)\\
{\phantom{1\over 2}}&-(x-\mu^i)^T\Sigma^i\nu^i R^i A^i x-x^T (A^i)^T R^iA^i x+\lambda^i=f^i(x;m^1,\ldots,m^N)
\end{align*}
or, equivalently,
\begin{align*}
-\Tr(\nu^iR^i\nu^i\Sigma^i+\nu^iR^iA^i)&+\,{1\over 2}\, (x-\mu^i)^T\Sigma^i\nu^i R^i\nu^i\Sigma^i(x-\mu^i)
-\,{1\over 2}\, x^T(A^i)^T R^iA^ix\\
&+\,{1\over 2}\Big(x^T  (A^i)^T R^i\nu^i\Sigma^i (x-\mu^i)-(x-\mu^i)^T \Sigma^i\nu^i R^i A^i x\Big)\\
{\phantom{1\over 2}}&+\lambda^i=f^i(x;m^1,\ldots,m^N)
\end{align*}
Taking into account that $R^i$, $\nu^i$, and $\Sigma^i$ are all symmetric matrices, the term in the second line vanishes and we obtain
\begin{align}\label{eq:quadratic}
-\Tr(\nu^iR^i\nu^i\Sigma^i+\nu^iR^iA^i)&+\, (x-\mu^i)^T\Sigma^i\,{\nu^i R^i\nu^i\over 2}\,\Sigma^i(x-\mu^i)\nonumber\\
&-\, x^T\,{(A^i)^T R^iA^i\over 2}\,x+\lambda^i=f^i(x;m^1,\ldots,m^N)
\end{align}


\n We can now exploit~\eqref{eq:ansatz_meas} to compute explicitly the expression of $f^i$. Indeed, since $m^i={\cal N}(\mu^i,(\Sigma^i)^{-1})$, we have
\begin{align*}
f^i(X^i;m^1,\ldots,m^N)&=\sum_{j,k=1}^N \int_{\R^{d(N-1)}} (X^j-\overline{X_i}^j)^T Q^i_{jk}(X^k-\overline{X_i}^k)\prod_{\ell\neq i} dm^\ell(X^\ell)\\
&= (X^i-\overline{X_i}^i)^T Q^i_{ii} (X^i-\overline{X_i}^i)+(X^i-\overline{X_i}^i)^T \Big(\sum_{k\neq i} Q^i_{ik}(\mu^k-\overline{X_i}^k)\Big)\\
&~~~~~~~~+\Big(\sum_{j\neq i} (\mu^j-\overline{X_i}^j)^T Q^i_{ji}\Big)(X^i-\overline{X_i}^i)\\
&~~~~~~~~+\sum_{j,k\neq i\,,j\neq k}(\mu^j-\overline{X_i}^j)^T Q^i_{jk}(\mu^k-\overline{X_i}^k)\\
&~~~~~~~~+\sum_{j\neq i}\left(\Tr(Q^i_{jj} (\Sigma^i)^{-1}) + (\mu^j-\overline{X_i}^j)^TQ^i_{jj}(\mu^j-\overline{X_i}^j)\right)\\
&=: (X^i)^T F^i_2 X^i+(X^i)^T F^i_{1,1}+F^i_{1,2} X^i +F^i_0
\end{align*}
where we have used the relation
$\vatt[v^TM v]=\Tr(M \Sigma^{-1})+(\mu)^TM\mu$,
which holds for any symmetric matrix $M$ and any vector of random variables $v$ whose expected value is $\mu$ and whose covariance matrix is $\Sigma^{-1}$,
to compute explicitly the last quadratic term. For later use, we write explicitly the expressions of $F^i_2$, $F^i_{1,1}$, $F^i_{1,2}$ and $F^i_0$:
$$
F^i_2= Q^i_{ii}\,,
$$
$$
F^i_{1,1}= -Q^i_{ii}\overline{X_i}^i+\Big(\sum_{j\neq i} Q^i_{ij}(\mu^j-\overline{X_i}^j)\Big)
$$
$$
F^i_{1,2}= -(\overline{X_i}^i)^TQ^i_{ii}+\Big(\sum_{j\neq i} (\mu^j-\overline{X_i}^j)^T Q^i_{ji}\Big)
$$
\begin{align*}
F^i_0&= (\overline{X_i}^i)^TQ^i_{ii}\overline{X_i}^i-(\overline{X_i}^i)^T\Big(\sum_{j\neq i} Q^i_{ij}(\mu^j-\overline{X_i}^j)\Big)
-\Big(\sum_{j\neq i} (\mu^j-\overline{X_i}^j)^T Q^i_{ji}\Big)\overline{X_i}^i\\
&~~~+\!\!\sum_{j,k\neq i\,,j\neq k}(\mu^j-\overline{X_i}^j)^T Q^i_{jk}(\mu^k-\overline{X_i}^k)
+\sum_{j\neq i}\left(\Tr(Q^i_{jj}  (\Sigma^i)^{-1}) + (\mu^j-\overline{X_i}^j)^TQ^i_{jj}(\mu^j-\overline{X_i}^j)\right)
\end{align*}


\n Once again, we can interpret equation~\eqref{eq:quadratic}, which is equivalent to the first equation of~\eqref{eq:hjkfp}, as an equality between quadratic forms to be satisfied for every $x\in\R^d$. This means that we must equate the corresponding coefficients. Notice that the assumption of $Q^i$ symmetric in~\eqref{eq:quad_cost} implies that $(Q^i_{jk})^T=Q^i_{kj}$ for every $j,k\in\{1,\ldots,N\}$. In particular, $Q^i_{ii}\in\Sym$. Hence, $F^i_{1,1}=(F^i_{1,2})^T$ and 
the two conditions on the linear terms
$$
-x^T\,{\Sigma^i \nu^iR^i\nu^i\Sigma^i\over 2}\, \mu^i=x^T F^i_{1,1}\,,
\qquad\qquad
-(\mu^i)^T\,{\Sigma^i \nu^iR^i\nu^i\Sigma^i\over 2}\, x=F^i_{1,2} x\,,
$$
do coincide. This 
leads to three conditions on the coefficients of the quadratic forms which have to be satisfied by the matrices $\Sigma^i$, the vectors $\mu^i$ and the real numbers $\lambda^i$.
\begin{align}\label{eq:cond1}
\Sigma^i\,{\nu^i R^i\nu^i\over 2}\,\Sigma^i-\,{(A^i)^T R^iA^i\over 2}\,=F^i_2
\end{align}
\begin{equation}\label{eq:cond2}
-\,{\Sigma^i \nu^iR^i\nu^i\Sigma^i\over 2}\, \mu^i=F^i_{1,1}
\end{equation}
\begin{equation}\label{eq:cond3}
(\mu^i)^T\,{\Sigma^i \nu^iR^i\nu^i\Sigma^i\over 2}\, \mu^i-\Tr(\nu^iR^i\nu^i\Sigma^i+\nu^iR^iA^i)+\lambda^i=F^i_0
\end{equation}


\n {\bf Step~3.} Notice that~\eqref{eq:cond1} is equivalent to say that $\Sigma^i$ solves an ARE of the form~\eqref{eq:ARE} with 
$$
{\cal R}=\,{\nu^i R^i\nu^i\over 2}\,,
\qquad\qquad
{\cal Q}=\,{(A^i)^T R^iA^i\over 2}\,+Q^i_{ii}\,.
$$ 
Hence, its solutions can be found as $d$--dimensional invariant graph subspaces of the $2d\times 2d$ matrix ${\cal H}$ defined in~\eqref{eq:inv_spaces}.
\n Since our standing assumptions {\bf (H)} imply that both ${\cal R}$ and ${\cal Q}$ are positive definite, we can apply Proposition~\ref{prop:ARE} {\it (ii)}--{\it (iii)} to conclude that~\eqref{eq:cond1} has a unique solution $\Sigma^i$ in $\Symp$, representing the inverse of the covariance matrix for our multivariate Gaussian $m^i$. 


\n {\bf Step~4.} Concerning~\eqref{eq:cond2}, we can rewrite the condition as
\begin{equation}\label{eq:eq2}
-\,{\Sigma^i \nu^iR^i\nu^i\Sigma^i\over 2}\,\mu^i-\sum_{j\neq i} Q^i_{ij} \mu^j=-\sum_{j=1}^N Q^i_{ij}\,\overline{X_i}^j 
\end{equation}
or equivalently, by collecting the relations~\eqref{eq:eq2} for $i=1,\ldots,N$,
\begin{equation}\label{eq:system_mi}
{\cal B}~\left(
\begin{array}{c}
\mu^1\\
\vdots\\
\mu^N
\end{array}
\right)= P
\end{equation}
where ${\cal B}$ is the $Nd\times Nd$ matrix defined in~\eqref{eq:matrixB} and $P$ is the vector in $\R^{Nd}$ defined by~\eqref{eq:vecP}. Here, we have also used the fact that each $\Sigma^p$ solves the Riccati equation~\eqref{eq:cond1} to rewrite the terms ${\cal B}_{pp}$ as
$$
{\cal B}_{pp}=-\,{\Sigma^p \nu^pR^p\nu^p\Sigma^p\over 2}=-\,{(A^p)^T R^pA^p\over 2}\,-Q^p_{pp}\,,
$$
i.e. in the form expected in~\eqref{eq:matrixB}.
%
%
\n Finally,~\eqref{eq:cond3} becomes
$$
\lambda^i=F^i_0-(\mu^i)^T\,{\Sigma^i \nu^iR^i\nu^i\Sigma^i\over 2}\, \mu^i+\Tr(\nu^iR^i\nu^i\Sigma^i+\nu^iR^iA^i)\,.
$$


\n {\bf Step~5.} So far we have mostly manipulated the equations of~\eqref{eq:hjkfp}, under the assumptions~\eqref{eq:quad_gauss}, arriving to an equivalent system
of matrix equations~\eqref{eq:KFP_matrix}, \eqref{eq:cond1}, \eqref{eq:cond2}, and~\eqref{eq:cond3}. 

 Now let us assume that condition {\bf (E)} holds. We have seen that each equation~\eqref{eq:cond1}, for $i\in\{1,\ldots,N\}$, is an ARE which admits a unique symmetric and positive definite solution $\Sigma^i$. By {\bf (E)}, the matrices $\Sigma^i$ also satisfy the Sylvester equation~\eqref{eq:sylvester} and thus
$$
R^i\nu^i\Sigma^i+R^iA^i=\Sigma^i\nu^iR^i+(A^i)^TR^i=\Big(R^i\nu^i\Sigma^i+R^iA^i\Big)^T\,,
$$
so that, by setting $\Lambda^i$ according to the first relation in~\eqref{eq:KFP_matrix}, we obtain a symmetric matrix as required. Moreover, the assumption on ${\cal B}$ and $[{\cal B}, P]$ in {\bf (E)} ensures the existence of a solution $(\mu^1,\ldots,\mu^N)$ to the linear system~\eqref{eq:system_mi} of $Nd$ equations in $Nd$ unknowns.
By using the solutions $\Sigma^i$ and $\mu^i$ in~\eqref{eq:cond3} and in the second relation of~\eqref{eq:KFP_matrix}, we also find admissible values $\lambda^i$ and $\rho^i$, and these complete the construction of a solution of the form~\eqref{eq:quad_gauss}.

 Viceversa, let us assume that a solution of the form~\eqref{eq:quad_gauss} exists for suitable matrices $\Lambda^i\in\Sym$, $\Sigma^i\in\Symp$, and suitable vectors $\rho^i,\mu^i$. Then, by the analysis in Step~1, we necessarily have~\eqref{eq:KFP_matrix}. Furthermore, by the analysis in Step~2, $\Sigma^i$ must solve the ARE~\eqref{eq:cond1} and $\mu^i$ and $\lambda^i$ must be given by solutions to~\eqref{eq:cond2} and~\eqref{eq:cond3}. In particular, the system~\eqref{eq:system_mi} admits at least a solution, and this implies the condition on the rank of ${\cal B}$ and $[{\cal B}, P]$. Finally, by combining the symmetry of $\Lambda^i$ with~\eqref{eq:KFP_matrix}, one has that the (unique) solution to~\eqref{eq:cond1} has to satisfy
$$
0=\Lambda^i-(\Lambda^i)^T=R^i\nu^i\Sigma^i+R^iA^i-\Big(R^i\nu^i\Sigma^i+R^iA^i\Big)^T\,,
$$
which is equivalent to~\eqref{eq:sylvester}. Thus, both requirements of {\bf (E)} must be necessarily satisfied.


\n {\bf Step~6.} We now focus our attention on the uniqueness. For~\eqref{eq:cond1} there is nothing to prove, because uniqueness of solution in $\Symp$ always follows from Proposition~\ref{prop:ARE}, under hypotheses {\bf (H)}. For~\eqref{eq:cond2}, the equivalence between {\bf (U)} and the uniqueness of solution $(\mu^1,\ldots,\mu^N)$ is evident when considering the equivalent form~\eqref{eq:system_mi}. Finally, once we have a unique choice for the matrices $(\Sigma^1,\ldots,\Sigma^N)$ and for the vectors $(\mu^1,\ldots,\mu^N)$, the uniqueness of $\lambda^i$, $\Lambda^i$ and $\rho^i$ is verified immediately.


\n {\bf Step~7.} It remains to prove that the affine feedbacks~\eqref{eq:feedback_strateg} provide a Nash equilibrium strategy for the game. Consider
$$
\overline{\alpha}^i(x)=(R^i)^{-1}\nabla v^i(x)=(R^i)^{-1}(\Lambda^ix+\rho^i) = (\nu^i \Sigma^i+A^i)\,x-\nu^i \Sigma^i\mu^i\,.
$$
By Proposition~\ref{prop:affine_admissible}, we know that $\overline{\alpha}^i(x)$ is admissible and that the ergodic measure associated to the process $X^i_t$ which solves
$$
dX^i_t=[-\nu^i \Sigma^iX^i_t+\nu^i \Sigma^i\mu^i]dt+\sigma^idW^i_t\,,
$$
is a multivariate Gaussian ${\cal N}(\mu,V)$ with mean
$$
\mu=(-\nu^i \Sigma^i)^{-1}(-\nu^i \Sigma^i\mu^i)=\mu^i\,,
$$
and covariance matrix $V=(\Sigma^i)^{-1}$, since
$$
(-\nu^i \Sigma^i)(\Sigma^i)^{-1}+(\Sigma^i)^{-1}(-\nu^i \Sigma^i)^T+\sigma^i(\sigma^i)^T=-2\nu^i+\sigma^i(\sigma^i)^T=0\,,
$$
and the equation~\eqref{eq:var_equation} admits a unique solution (see Remark~\ref{rem:var_equation}). In other words, the invariant measure coincides with the measure $m^i$ satisfying~\eqref{eq:ansatz_meas}.

 We now consider any other admissible strategy $\alpha^i$ and obtain from Dynkin--It\^o's formula
$$
\vatt\left[v^i(X^i_T)-v^i(X^i_0)\right]=
\vatt\left[\int_0^T\nabla v^i(X^i_s)\cdot(A^iX^i_s-\alpha^i_s)\,+\Tr\left({\sigma^i(\sigma^i)^T\over 2}\,\D^2 v^i(X^i_s)\right)\,ds\right]\,.
$$
Hence, from $\nu^i=\sigma^i(\sigma^i)^T/2$ and the fact that the map $(x,y)\mapsto x^TR^iy$ is an inner product, one obtains
\begin{align*}
\vatt\left[v^i(X^i_T)-v^i(X^i_0)\right]&=
\vatt\left[\int_0^T\Big(\Tr(\nu^i\D^2 v^i)+(\nabla v^i)^TA^ix-(\nabla v^i)^T\alpha^i_s\Big)(X^i_s)\,ds\right]\\
&=\vatt\left[\int_0^T\Big(\Tr(\nu^i\D^2 v^i)+(\nabla v^i)^TA^ix-\big((R^i)^{-1}\nabla v^i\big)^TR^i\alpha^i_s\Big)(X^i_s)\,ds\right]\\
&\geq\vatt\left[\int_0^T\Big(\Tr(\nu^i\D^2 v^i)+(\nabla v^i)^TA^ix-\,{1\over 2}\,(\nabla v^i)^T(R^i)^{-1}\nabla v^i\Big)(X^i_s)-\,{1\over 2}\,(\alpha^i_s)^T R^i\alpha^i_s\,ds\right]\\
&=\vatt\left[\int_0^T\Tr(\nu^i\D^2 v^i(X^i_s))-H^i(X^i_s,\nabla v^i(X^i_s))-\,{1\over 2}\,(\alpha^i_s)^T R^i\alpha^i_s\,ds\right]
\end{align*}
with equality holding if $\alpha^i=\overline{\alpha}^i$. Therefore, the first equation in~\eqref{eq:hjkfp} implies
\begin{align*}
\vatt\left[v^i(X^i_T)-v^i(X^i_0)\right]&\geq
\vatt\left[\int_0^T\lambda^i -f^i(X^i_s)-\,{1\over 2}\,(\alpha^i_s)^T R^i\alpha^i_s\,ds\right]\\
&=\lambda^iT- \vatt\left[\int_0^Tf^i(X^i_s)+\,{1\over 2}\,(\alpha^i_s)^T R^i\alpha^i_s\,ds\right]
\end{align*}
Hence, by dividing by $T$ and letting $T\to+\infty$, we get
\begin{equation}\label{lambda_i}
\lambda^i\leq\liminf_{T\to+\infty}\,{1\over T}~\vatt\left[\int_0^Tf^i(X^i_s)+\,{1\over 2}\,(\alpha^i_s)^T R^i\alpha^i_s\,ds\right]\,,
\end{equation}
because the left hand side of the original inequality vanishes due to $v^i$ being quadratic and the strategies being admissible (and therefore $\vatt[X^i_t]\leq C$ and 
$\vatt[(X^i_t)(X^i_t)^T]\leq C$
 for some constant $C$). 

It 
remains to prove that the right hand side of \eqref{lambda_i} is $J^i(X,\overline{\alpha}^1,\ldots\overline{\alpha}^{i-1},\alpha^i,\overline{\alpha}^{i+1},\ldots,\overline{\alpha}^N)$, 
which means that the cost $\lambda^i$ corresponds to a Nash equilibrium.
%
 This 
  property 
  follows from the ergodicity. Indeed, let us consider the probability measures
$$
\tilde m^i:= m^i=m^i(\alpha^i)\,,
\qquad\qquad
\tilde m^j:= m^j=m^j(\overline{\alpha}^j) \qquad\forall~j\neq i\,,
$$
which are the invariant measures corresponding to the solution $(X^1_t,\ldots,X^N_t)$ obtained when players adopt the strategies $(\overline{\alpha}^1,\ldots\overline{\alpha}^{i-1},\alpha^i,\overline{\alpha}^{i+1},\ldots,\overline{\alpha}^N)$. By recalling that the functions $f^i$ are polynomials of degree less or equal to two in the variable $x^i$, we obtain
\begin{align*}
{1\over T}~\vatt\left[\int_0^T f^i(X^i_s)\,ds\right]&={1\over T}~\vatt\left[\int_0^T\left(\int_{\R^{d(N-1)}} (X-\overline{X_i})^T Q^i(X-\overline{X_i})\prod_{\ell\neq i} d\tilde m^\ell(X^\ell)\right)(X^i_s)\,ds\right]\\
&={1\over T}~\vatt\left[\int_0^T h(X^i_s)\,ds\right]\xrightarrow{T\to\infty}\int_{\R^d}h(x)\,d\tilde m^i(x)\\
&=\int_{\R^{Nd}} \Big((X-\overline{X_i})^T Q^i(X-\overline{X_i})\Big)\prod_{\ell=1}^N d\tilde m^\ell(X^\ell)\\
&=\sum_{j=1}^N\int_{\R^{d}} \Big( (X^j-\overline{X_i}^j)^T Q^i_{jj}(X^j-\overline{X_i}^j)\Big) \,d\tilde m^j(X^j)\\
&~~~~~~~+\sum_{j\neq k}\int_{\R^{2d}} \Big( (X^j-\overline{X_i}^j)^T Q^i_{jk}(X^k-\overline{X_i}^k)\Big)\, d\tilde m^j(X^j)\,d\tilde m^k(X^k)\\
&=\lim_{T\to+\infty}\,{1\over T}~\vatt \left[\sum_{j=1}^N\int_0^T \Big( (X^j_s-\overline{X_i}^j)^T Q^i_{jj}(X^j_s-\overline{X_i}^j)\Big) \,ds\right]\\
&~~~~~~~+\lim_{T\to+\infty}\,{1\over T}~\vatt \left[\sum_{j\neq k}\int_0^T \Big( (X^j_s-\overline{X_i}^j)^T Q^i_{jk}(X^k_s-\overline{X_i}^k)\Big)\, ds\right]\\
&=\lim_{T\to+\infty}\,{1\over T}~\vatt \left[\int_0^T F^i(X^1_s,\ldots,X^N_s)\,ds\right]
\end{align*}
where we have also used the ergodicity of the pair $(X^j_t,X^k_t)$ with corresponding measure given by the product measure obtained from $\tilde m^j$ and $\tilde m^k$. In conclusion, 
the right hand side of \eqref{lambda_i} is $J^i(X,\overline{\alpha}^1,\ldots\overline{\alpha}^{i-1},\alpha^i,\overline{\alpha}^{i+1},\ldots,\overline{\alpha}^N)$ and this completes the proof.~~$\diamond$

\begin{rem} \rm{
Looking at the formulas~\eqref{eq:quad_gauss} and~\eqref{eq:KFP_matrix} it might seem that the Nash equilibrium strategies depend on the noise $\sigma^i$, through $\nu^i$, which is typically not the case in LQ stochastic problems.
In fact, this is not the case in our problem either: by introducing new variables ${\cal V}^i:=\nu^i\Sigma^i$, it is immediate to verify that the feedback strategies only depend on these new variables ${\cal V}^i$, which are determined by the Riccati equations
$$
({\cal V}^i)^T\,{R^i\over 2}\,{\cal V}^i=Q_{ii}^i+\,{(A^i)^T R^iA^i\over 2}\,,
$$
and thus do not depend on the noise statistics $\sigma^i$. 
The same holds for 
the 
 equilibrium strategies 
of  Theorem~\ref{thm:NIDplay} in the next Section. This allows to take the small noise or vanishing viscosity limit $\nu^i\to0$ under some additional conditions, extending the results on the case $d=1$ in \cite{Bardi}, see the sequel of this paper \cite{Priuli}.
}
\end{rem}

\section[Nearly identical players]{Nearly identical players} \label{sec:nearly_id}

In this section we introduce assumptions saying that the players are almost identical, as in~\cite{Bardi}, and prove that there exists a Nash equilibrium with the same feedback and the same distribution for all players, although the values can be different.
The first condition is a \emph{Symmetry assumption} on the cost of each player:

\begin{description}
\item{\bf (S)} every player is influenced in the same way by other players, i.e. for each $i\in\{1,\ldots,N\}$ and each $j,k\neq i$
$$
F^i(X^1,\ldots,X^j,\ldots, X^k,\ldots, X^N)=F^i(X^1,\ldots,X^k,\ldots, X^j,\ldots, X^N)
$$
\end{description}

\n We can easily prove the following lemma.

\begin{lem}\label{lem:symmetry} Assumption {\bf (S)} holds if and only if there exist matrices $B_i,C_i,D_i$ and vectors $\Delta_i$ such that
$$
Q^i_{ij}=\,{B_i\over 2}\,,
\qquad
Q^i_{jj}=C_i\,,
\qquad
\overline{X_i}^j=\Delta_i\,,
\qquad
\qquad
\forall~j\neq i\,,
$$
$$
Q^i_{jk}=D_i\,,
\qquad
\qquad
\forall~j,k\neq i\,,j\neq k\,.
$$
\end{lem}

\n Under assumption {\bf (S)}, the quadratic costs $F^i$ take the following form
\begin{align*}
F^i(X^1,\ldots, X^N)&= (X^i-\overline{X_i}^i)^T Q^i_{ii} (X^i-\overline{X_i}^i)\\
&~~~~~+(X^i-\overline{X_i}^i)^T \,{B_i\over 2}\,\Big(\sum_{k\neq i} (X^k-\Delta_i)\Big)+\Big(\sum_{j\neq i} (X^j-\Delta_i)^T\Big)\,{B_i\over 2}\,(X^i-\overline{X_i}^i)\\
&~~~~~+\sum_{j\neq i}(X^j-\Delta_i)^T C_i(X^j-\Delta_i)+\sum_{j,k\neq i\,,j\neq k}(X^j-\Delta_i)^T D_i(X^k-\Delta_i)
\end{align*}
In particular, they can be written in the form usually arising in mean field games, namely
$$
F^i(X^1,\ldots, X^N)=V^i\left[{1\over N-1}\sum_{j\neq i}\delta_{X^j}\right](X^i)\,,
$$
where $\delta_{X^j}$ is the Dirac measure on $\R^d$ centered in the point $X^j$ and $V^i$ is the operator, mapping probability measures ${\frak m}$ on $\R^d$, with finite second moments, into quadratic polynomials, defined by the expression
\begin{align*}
V^i[{\frak m}](X)&:= (X^i-\overline{X_i}^i)^T Q^i_{ii} (X^i-\overline{X_i}^i)\\
&~~~~~+(N-1)\int_{\R^d}\left((X^i-\overline{X_i}^i)^T \,{B_i\over 2}\, (\xi-\Delta_i)+(\xi-\Delta_i)^T\,{B_i\over 2}\,(X^i-\overline{X_i}^i)\right)\,d{\frak m}(\xi)\\
&~~~~~+(N-1)\int_{\R^d} (\xi-\Delta_i)^T (C_i-D_i) (\xi-\Delta_i)\,d{\frak m}(\xi)\\
&~~~~~+\left((N-1)\int_{\R^d} (\xi-\Delta_i)\,d{\frak m}(\xi)\right)^TD_i\left((N-1)\int_{\R^d} (\xi-\Delta_i)\,d{\frak m}(\xi)\right)
\end{align*}
Indeed, it is enough to recall that for any choice of vectors $w_1,\ldots,w_N$ and of an index $i\in\{1,\ldots,N\}$, there holds
$$
\sum_{j,k\neq i\,,j\neq k}w_j^T D_iw_k=\left(\sum_{j\neq i} w_j^T\right) D_i\left(\sum_{j\neq i} w_j\right)-\sum_{j\neq i} w_j^TD_i w_j\,.
$$

\begin{defn}\label{def:nearly_id} We say that the players are \emph{nearly identical} if costs $F^i$ satisfy {\bf (S)} and if all players have the same:
\begin{itemize}
\item control systems, i.e. $A^i=A$ and $\sigma^i=\sigma$ (and therefore $\nu^i=\nu$) for all $i$,
\item costs of the control, i.e. $R^i=R$ for all $i$,
\item reference positions, i.e. $\overline{X_i}^i=H$ (own reference position, or happy place) and $\Delta_i=\Delta$ (reference position of the other players) for all $i$,
\item primary costs of displacement, i.e. $Q^i_{ii}=Q$ and $B_i=B$ for all $i$.
\end{itemize}
\end{defn}

Note that the players are not fully identical because the secondary costs of displacement 
 $C_i$ and $D_i$ can 
  be different among them. 
   Observe that in this framework, the hypotheses {\bf (H)} specialize to
\begin{equation}\label{eq:NIDplay_hyp}
\det(\sigma)\neq 0\,,
\qquad\qquad
R\in\Symp\,,
\qquad\qquad
Q\in\Symp\,.
\end{equation}

 Let us rewrite part of the computations from the previous section for nearly identical players. First of all, the 
 right hand side $f^i$ \eqref{eq:rhs} of the HJB equation becomes, for given measures $m^i={\cal N}(\mu^i,(\Sigma^i)^{-1})$, 
\begin{align*}
f^i(X^i;m^1,\ldots,m^N)&= (X^i-H)^T Q (X^i-H)\\
&~~~+(X^i-H)^T\,{B\over 2}\, \left(\sum_{k\neq i} (\mu^k-\Delta)\right)+\left(\sum_{j\neq i} (\mu^j-\Delta)^T\right)\,{B\over 2}\,(X^i-H)\\
&~~~+(N-1)~\Tr((C_i-D_i)\,  (\Sigma^i)^{-1}) +\sum_{j\neq i} (\mu^j-\Delta)^T(C_i-D_i)(\mu^j-\Delta)\\
&~~~+\left(\sum_{j\neq i}(\mu^j-\Delta)^T\right) D_i\left(\sum_{j\neq i}(\mu^j-\Delta)\right)
\end{align*}
Hence, if we search for identically distributed solutions for all 
 players, i.e. if we search for measures of the form $m^1=\ldots=m^N={\cal N}(\mu,\Sigma^{-1})$, we obtain
\begin{align*}
f^i(X^i;m^1,\ldots,m^N)&= (X^i-H)^T Q (X^i-H)\\
&~~~+(N-1)~(X^i-H)^T\,{B\over 2}\, (\mu-\Delta)+(N-1)~(\mu-\Delta)^T\,{B\over 2}\,(X^i-H)\\
&~~~+(N-1)~\Tr((C_i-D_i)\, \Sigma^{-1}) +(N-1)~ (\mu-\Delta)^TC_i(\mu-\Delta)\\
&~~~+(N-1)(N-2)~(\mu-\Delta)^T D_i(\mu-\Delta)
\end{align*}
Let us investigate the existence of solutions such that
\begin{equation}\label{eq:ansatz_nearly_id}
v^i(x) =x^T\,{ \Lambda\over 2}\, x+\rho^Tx\,,
\qquad
m^i(x)=\gamma\exp\left\{-\,{1\over 2} (x-\mu)^T\Sigma(x-\mu)\right\}
\end{equation}
for suitable symmetric matrices $\Lambda,\Sigma$, with $\Sigma$ positive definite, suitable vectors $\mu,\rho$ and a suitable constant $\gamma$ depending only on the matrix $\Sigma$ and on the dimension of the space. By repeating the same computations done in Section~\ref{sec:n_play}, it is immediate to verify that the KFP equation in \eqref{eq:hjkfp} for the measure reduces, as 
in \eqref{eq:KFP_matrix}, to the matrix relations
\begin{equation}\label{eq:KFP_matrix_nearly_id}
\Lambda=R\big(\nu\Sigma+A\big)\,,\qquad\qquad \rho=-R\nu\Sigma\mu\,,
\end{equation}
By interpreting again the HJB equation for the value function  in \eqref{eq:hjkfp} as an equality between quadratic forms 
we obtain the system of equations 
\begin{equation}\label{eq:ARE_nearly_id}
\Sigma\,{\nu R\nu\over 2}\,\Sigma-\,{A^T RA\over 2}\,=Q
\end{equation}
\begin{equation}\label{eq:mu_nearly_id}
\,-{\Sigma \nu R\nu \Sigma\over 2}\, \mu=-Q H+(N-1)~{B\over 2}\,(\mu-\Delta)
\end{equation}
\begin{equation}\label{eq:lambda_nearly_id}
\mu^T\,{\Sigma \nu R \nu \Sigma\over 2}\, \mu-\Tr(\nu R\nu\Sigma+\nu RA)+\lambda^i=\widetilde F^i_0
\end{equation}
with
\begin{align*}
\widetilde F^i_0&= H^TQ H-(N-1)~\left(H^T{B\over 2}\,(\mu-\Delta)+(\mu-\Delta)^T {B\over 2}\,H\right)\\
&~~~+(N-1)~\Tr((C_i-D_i) \Sigma^{-1}) +(N-1)~ (\mu-\Delta)^TC_i(\mu-\Delta)\\
&~~~+(N-1)(N-2)~(\mu-\Delta)^T D_i(\mu-\Delta)
\end{align*}

In particular, the first equation \eqref{eq:ARE_nearly_id} has exactly the same form as~\eqref{eq:cond1}, i.e. it is again an ARE which admits a unique solution $\Sigma$ in $\Symp$, under hypotheses~\eqref{eq:NIDplay_hyp}. By plugging \eqref{eq:ARE_nearly_id} 
into \eqref{eq:mu_nearly_id} this can be rewritten as 
$$
-\left(Q+\,{A^TRA\over 2}\,-(1-N)~{B\over 2}\right)\, \mu=-Q H+(1-N)~{B\over 2}\,\Delta\,,
$$
which admits a unique solution $\mu\in\R^d$ whenever the matrix
\begin{equation}\label{eq:weaker_invert}
{\cal B}':= Q+\,{A^TRA\over 2}\,-(1-N)~{B\over 2}
\end{equation}
is invertible. Finally, once $\Sigma$ and $\mu$ have been found, they can be used in the third equation~\eqref{eq:lambda_nearly_id} and in~\eqref{eq:KFP_matrix_nearly_id} to obtain the values $\lambda^i\in\R,\, i=1,\dots, N$, the matrix $\Lambda$ and the vector $\rho$. 

For nearly identical players  the appropriate analogs of 
the conditions 
{\bf (E)} and {\bf (U)} are the following.

\begin{description}
\item{\bf (E$'$)} Every symmetric and positive definite solution $Y$ of the algebraic Riccati equation
\begin{equation}\label{eq:riccati_nearlyid}
Y\,{\nu R\nu\over 2}\,Y=\,{A^T RA\over 2}\,+Q\,,
\end{equation}
is also a solution to the Sylvester equation
\begin{equation}\label{eq:sylvester_nearlyid}
Y \nu R-R\nu Y= R A-A^TR\,.
\end{equation}
Moreover, the matrices ${\cal B}'\in\Mat_{d\times d}(\R)$ and $[{\cal B}', P']\in\Mat_{d\times (d+1)}(\R)$ have the same rank, where ${\cal B}'$ is the matrix defined in~\eqref{eq:weaker_invert}, $P':=-Q H+(1-N)~{B\over 2}\,\Delta$, and $[{\cal B}', P']$ is defined as in
~\eqref{eq:expanded_matrix}.
\end{description}

\begin{description}
\item{\bf (U$'$)} The matrix ${\cal B}'$ defined in~\eqref{eq:weaker_invert} is invertible.
\end{description}

It is immediate to verify that the analysis performed above proves 
the following theorem. 

\begin{thm}\label{thm:NIDplay} Consider an $N$--players game with dynamics~\eqref{eq:sde} and costs~\eqref{eq:cost}. Assume that players are nearly identical and that~\eqref{eq:NIDplay_hyp} holds. Then, the associated system of $2N$ HJB--KFP equations~\eqref{eq:hjkfp} admits solutions $(v^i,m^i,\lambda^i)$, $i=1,\ldots,N$, 
with $v^1=\ldots=v^N=v$ quadratic function and $m^1=\ldots=m^N={\cal N}(\mu,\Sigma^{-1})$ multivariate Gaussian, i.e. of the form~\eqref{eq:ansatz_nearly_id}, if and only if condition {\bf (E$'$)} is satisfied. 
{In particular, $\Sigma$ is the unique solution in $\Symp$ of the Riccati equation \eqref{eq:riccati_nearlyid}, $\mu$ solves \eqref{eq:mu_nearly_id}, 
$\Lambda$ and $\rho$ are given by \eqref{eq:KFP_matrix_nearly_id}, and $\lambda^i$ by \eqref{eq:lambda_nearly_id}.}

\n Moreover, solutions of the form~\eqref{eq:ansatz_nearly_id} are also unique if and only if condition {\bf (U$'$)} is satisfied and, if this is the case, the affine feedbacks
$$
\overline{\alpha}^i(x)=\overline{\alpha}(x):= R^{-1}\nabla v(x)\,,\qquad\qquad x\in\R^d,~~i=1,\ldots,N
$$
provide a Nash equilibrium strategy for all initial positions $X\in\R^{Nd}$, among the admissible strategies, and $J^i(X,\overline{\alpha})=\lambda^i$ for all $X$ and all $i=1,\dots, N$.
\end{thm}
\begin{rem} \rm{Note that the distribution $m$ and the solution $v$ found in this theorem are the same even if the cost functionals of the players differ in the terms involving the matrices $C_i$ and $D_i$. These terms only affect the values of the game $\lambda^i$ 
{(this motivates the name ``secondary costs'' given to them). This result supports the existence of a large population limit, that we study in the next section.}
}\end{rem}
\begin{rem} \rm{
Note that the invertibility of ${\cal B}$ in~\eqref{eq:matrixB} is a stronger requirement than the invertibility of the matrix ${\cal B}'$ in~\eqref{eq:weaker_invert}. Indeed, if we take a game such that
$
Q+{A^TRA
/ 2}={B
/ 2}\,,
$ 
then ${\cal B}$ consists of blocks ${\cal B}_{\alpha\beta}=-B\in\Mat_{d\times d}(\R)$ for all $\alpha,\beta=1,\ldots,N$, and thus is not invertible since it satisfies $\mathrm{rank}({\cal B})=\mathrm{rank}(B)\leq d$. On the other hand, the matrix 
${\cal B}'=NB/2$ is invertible, provided $B$ is invertible. Therefore there can be infinitely many quadratic-Gaussian solutions although only one of them is identically distributed, see \cite{Bardi} for an explicit example.
}\end{rem}


\section[The large population limit
]{The large population limit
} \label{sec:limit}

In this section we study the convergence of Nash equilibria when the number $N$ of players goes to infinity. 
Assume for simplicity that the control system, the costs of the control and the reference positions are always the same, i.e. that $A,\sigma,R,H$ and $\Delta$ are all independent from the number of players $N$. We denote with
$$
Q^N\,,
\qquad
B^N\,,
\qquad
C_i^N\,,
\qquad
D_i^N\,,
$$
the primary and secondary costs of displacement, respectively, which are assumed to depend on $N$. We assume that these quantities, when $N\to+\infty$, tends to suitable matrices $\hat Q,\hat B,\hat C,\hat D$ with their natural scaling, i.e. as $N\to+\infty$ there holds
\begin{equation}\label{eq:scale}
Q^N\to \hat Q\,,
\qquad
B^N(N-1)\to \hat B\,,
\qquad
C_i^N(N-1)\to \hat C\,,
\qquad
D_i^N(N-1)^2\to \hat D\,,
\qquad
\forall~i\,.
\end{equation}
We define an operator acting on probability measures with finite second moments 
${\frak m}\in{\cal P}_2(\R^d)$ that describes the cost for an average player of the density ${\frak m}$ of the other players
\begin{align*}
\hat V[{\frak m}](X)&:=(X-H)^T \hat Q (X-H)\\
&~~~~~+\int_{\R^d}\left((X-H)^T \,{\hat B\over 2}\, (\xi-\Delta)+(\xi-\Delta)^T\,{\hat B\over 2}\,(X-H)\right)\,d{\frak m}(\xi)\\
&~~~~~+\int_{\R^d} (\xi-\Delta)^T \hat C (\xi-\Delta)\,d{\frak m}(\xi)\\
&~~~~~+\left(\int_{\R^d} (\xi-\Delta)\,d{\frak m}(\xi)\right)^T\hat D\left(\int_{\R^d} (\xi-\Delta)\,d{\frak m}(\xi)\right).
\end{align*}
Since both $V^i_N[{\frak m}]$ and $\hat V[{\frak m}]$ are quadratic forms on $\R^d$, for all $i$ and all ${\frak m}\in{\cal P}_2(\R^d)$, it is immediate to deduce that the convergence of the matrix coefficients in~\eqref{eq:scale} implies, as $N\to+\infty$,
$$
V^i_N[{\frak m}](X) \to \hat V[{\frak m}](X)\,,~~~\mbox{locally uniformly in $X$}\,.
$$


\n By denoting with $\lambda^i_N$, $v_N$ and $m_N$ the solutions found in Theorem~\ref{thm:NIDplay} of Section~\ref{sec:nearly_id}, we expect that the limits of these solutions satisfy, like in~\cite{Bardi, LL1, LL3}, the system of two mean field HJB-KFP equations
\begin{equation}\label{eq:hjkfp_limit}
\left\{
\begin{array}{l}
-\Tr(\nu\D^2 v)+\displaystyle\,{1\over 2}\,\nabla v^T R^{-1}\nabla v-\nabla v^TAx+\lambda=\hat V[m](x)\\
~\\
-\Tr(\nu\D^2 m)-\displaystyle\Div\Big(m \cdot(R^{-1}\nabla v-A x)\Big)=0\\
~\\
\int_{\R^d}m(x)\,dx=1\,,\qquad m>0
\end{array}
\right.
\end{equation}
Along the lines of hypotheses {\bf (H)}, the natural assumptions on the coefficients of~\eqref{eq:hjkfp_limit} are the following
\begin{equation}\label{eq:hyp_limit}
\nu\in\Symp\,,
\qquad\qquad
R\in\Symp\,,
\qquad\qquad
\hat Q\in\Symp\,.
\end{equation}

\smallskip

We look for solutions of~\eqref{eq:hjkfp_limit} such that
\begin{equation}\label{eq:ansatz_limit}
v(x) =x^T\,{ \Lambda\over 2}\, x+\rho^T x\,,
\qquad\qquad
m(x)=\gamma\exp\left\{-\,{1\over 2} (x-\mu)^T\Sigma(x-\mu)\right\}
\end{equation}
for suitable symmetric matrices $\Lambda,\Sigma$, with $\Sigma$ positive definite, suitable vectors $\mu,\rho$ and a suitable constant $\gamma$ depending only on the matrix $\Sigma$ and on the dimension of the space. 

\n By repeating the computations done in Sections~\ref{sec:n_play} and~\ref{sec:nearly_id}, it is immediate to verify that the KFP equation  for the measure in~\eqref{eq:hjkfp_limit} reduces 
to the matrix equations 
\begin{equation}\label{eq:KFP_matrix_limit}
\Lambda=R\big(\nu\Sigma+A\big)\,,\qquad\qquad \rho=-R\nu\Sigma\mu\,.
\end{equation}
Concerning the HJB equation for the value function, one can proceed as in the previous sections, obtaining the system of matrix equations 
\begin{equation}\label{eq:ARE_limit}
\Sigma\,{\nu R\nu\over 2}\,\Sigma-\,{A^T RA\over 2}\,=\hat Q
\end{equation}
\begin{equation}\label{eq:mu_limit}
-\,{\Sigma \nu R\nu \Sigma\over 2}\, \mu=-\hat Q H+\,{\hat B\over 2}\,(\mu-\Delta)
\end{equation}
\begin{equation}\label{eq:lambda_limit}
\mu^T\,{\Sigma \nu R \nu \Sigma\over 2}\, \mu-\Tr(\nu R\nu\Sigma+\nu RA)+\lambda=\hat F_0
\end{equation}
with
\begin{align*}
\hat F_0&= H^T\hat Q H-\left(H^T{\hat B\over 2}\,(\mu-\Delta)+(\mu-\Delta)^T {\hat B\over 2}\,H\right)
+\Tr(\hat C \Sigma^{-1}) +(\mu-\Delta)^T(\hat C+\hat D)(\mu-\Delta)
\end{align*}

In particular, under assumptions~\eqref{eq:hyp_limit}, the first equation is an ARE which admits a unique solution $\Sigma$ in $\Symp$. Also, we can rewrite the second equality in the form
$$
-\left(\hat Q+\,{A^TRA\over 2}\,+\,{\hat B\over 2}\right)\, \mu=-\hat Q H-\,{\hat B\over 2}\,\Delta\,,
$$
which admits a unique solution $\mu$ whenever the matrix
\begin{equation}\label{eq:limit_invert}
{\cal B}_\infty:=\hat Q+\,{A^TRA\over 2}\,+\,{\hat B\over 2}
\end{equation}
is invertible. Finally, once $\Sigma$ and $\mu$ have been found, one can insert them into the third equation and~\eqref{eq:KFP_matrix_limit} to obtain the value $\lambda$, the matrix $\Lambda$ and the vector $\rho$ required by~\eqref{eq:ansatz_limit}.

\n In this case, the existence and uniquess of solutions to~\eqref{eq:hjkfp_limit} is then related to the following conditions.

\begin{description}
\item{\bf ($\mbox{E}_\infty$)} 
The symmetric and positive definite solution $Y$ of the algebraic Riccati equation
\begin{equation}\label{eq:riccati_limit}
Y\,{\nu R\nu\over 2}\,Y=\,{A^T RA\over 2}\,+\hat Q\,,
\end{equation}
is also a solution to the Sylvester equation
\begin{equation}\label{eq:sylvester_limit}
Y \nu R-R\nu Y= R A-A^TR\,.
\end{equation}
Moreover, the matrices ${\cal B}_\infty\in\Mat_{d\times d}(\R)$ and $[{\cal B}_\infty, P_\infty]\in\Mat_{d\times (d+1)}(\R)$ have the same rank, where ${\cal B}_\infty$ is the matrix defined in~\eqref{eq:limit_invert}, $P_\infty:=-\hat Q H+\,{\hat B\over 2}\,\Delta\in\R^d$ and $[{\cal B}_\infty, P_\infty]$ is defined analogously to~\eqref{eq:expanded_matrix}.
\end{description}

\begin{description}
\item{\bf ($\mbox{U}_\infty$)} The matrix ${\cal B}_\infty$ defined in~\eqref{eq:limit_invert} is invertible.
\end{description}

The main results of this section are the following two theorems.

\begin{thm}\label{thm:LIMplay} 
\begin{description}
\item{\it (i)}~{\bf [Solutions to MFPDE]} Assume
~\eqref{eq:hyp_limit}. 
 Then, the system of HJB--KFP equations~\eqref{eq:hjkfp_limit} admits solutions $(v,m,\lambda)$ with 
  $v$ quadratic function and $m$ multivariate Gaussians ${\cal N}(\mu,\Sigma^{-1})$, i.e. of the form~\eqref{eq:ansatz_limit}, if and only if condition {\bf ($\mbox{E}_\infty$)} is satisfied. In particular, $\Sigma$ is the unique solution in $\Symp$ of the Riccati equation \eqref{eq:ARE_limit}, $\mu$ solves \eqref{eq:mu_limit}, 
$\Lambda$ and $\rho$ are given by \eqref{eq:KFP_matrix_limit} and $\lambda$ by \eqref{eq:lambda_limit}.

\n Moreover, solutions of the form~\eqref{eq:ansatz_limit} are also unique if and only if condition {\bf ($\mbox{U}_\infty$)} is satisfied.

\item{\it (ii)}~{\bf [Convergence as $N\to\infty$]} Consider a sequence of differential games of the form~\eqref{eq:sde}--\eqref{eq:cost} with $N$ nearly identical players 
 and assume that~\eqref{eq:scale} is verified as $N\to\infty$. Also, assume that
 ~\eqref{eq:NIDplay_hyp} and {\bf ($\mathbf{E'
 }$)} with $Q=Q^N$ hold for all $N\in\N$, and that 
 ~\eqref{eq:hyp_limit}, {\bf ($\mbox{E}_\infty$)}, and {\bf ($\mbox{U}_\infty$)} are satisfied. Then, the solutions $(v_N,m_N,\lambda_N^1,\ldots,\lambda_N^N)$ found in Theorem \ref{thm:NIDplay} converge as $N\to\infty$
to the quadratic-Gaussian solution $(v,m,\lambda)$ of the mean-field system~\eqref{eq:hjkfp_limit} found in {\it (i)} 
 in the following sense: $v_N\to v$ in $C^1_{loc}(\R^d)$ with second derivative converging uniformly in $\R^d$, $m_N\to m$ in $C^k(\R^d)$ for all $k$, and $\lambda_N^i\to \lambda$ for all $i$.
\end{description}
\end{thm}

The uniqueness statement in 
 the previous theorem is only among solutions of \eqref{eq:hjkfp_limit} with $v$  quadratic and $m$ Gaussian. A natural question is whether the HJB-KFP system of PDEs 
  admits other $C^2$ solutions that are not of the form \eqref{eq:ansatz_limit} 
   and therefore other mean-field equilibria.
   We add a normalization condition on $v$, to avoid addition of constants, and make a simple assumption that ensures the monotonicity of $\hat V$ with respect to the scalar product in the Lebesgue space $L^2$. Then an argument of Lasry and Lions 
   ~\cite{LL1, LL3} implies uniqueness among general solutions satisfying natural growth conditions at infinity.


\begin{thm}\label{thm:uniq} 
The integral operator $\hat V$ satisfies $\int_{\R^d} \left(\hat V[{\frak m}]-\hat V[{\frak n}]\right)(x)\, d({\frak m}-{\frak n})(x)\geq 0$ for all probability measures ${\frak m},{\frak n}\in{\cal P}_2(\R^d)$ if and only if the matrix $\hat B$ is positive semidefinite.

\n Then 
for $\hat B\geq 0$ there is at most one solution $(v,m,\lambda)$  of~\eqref{eq:hjkfp_limit} with  $m\in{\cal P}_2(\R^d)$, $|\nabla v(x)|\leq C(1+|x|)$ for some $C>0$, and such that  $v(0)=0$;  in particular, under the assumptions of Theorem~\ref{thm:LIMplay},   the solution  
given by that theorem is the unique solution with such properties.
\end{thm}

\begin{rem} \rm{
The condition $B^N>0$ means that in the $N$-person game there is a positive cost for the $i$-th player if his displacement with respect to his happy state $H$ is the same as the average displacement of the other players from their reference position $\Delta$. Therefore the condition $\hat B\geq 0$ in the last theorem means that imitation among players is not rewarding in the large population limit. Note also that $\hat B\geq 0$ implies {\bf ($\mbox{U}_\infty$)} by \eqref{eq:hyp_limit}, but not viceversa. 
}\end{rem}

\n{\it Proof of Theorem~\ref{thm:LIMplay}.} The analysis performed on equations~\eqref{eq:KFP_matrix_limit}--\eqref{eq:lambda_limit} already proves part {\it (i)}. 
It remains to prove part {\it (ii)}. 
Let $(v_N,m_N,\lambda_N^1,\ldots,\lambda_N^N)$ be a solution of the differential game for $N$ nearly identical players found in Theorem~\ref{thm:NIDplay}. Then, 
\begin{equation}\label{formulaN}
v_N=x^T\,{ \Lambda_N\over 2}\, x+(\rho_N)^T x\,,
\qquad\qquad
m_N(x)=\gamma_N\exp\left\{-\,{1\over 2} (x-\mu_N)^T\Sigma_N(x-\mu_N)\right\}\,,
\end{equation}
where $\Sigma_N$ and $\mu_N$ solve equations~\eqref{eq:ARE_nearly_id} and \eqref{eq:mu_nearly_id}, respectively, $\Lambda_N$ and $\rho_N$ are given in terms of $\Sigma_N$ and $\mu_N$ by~\eqref{eq:KFP_matrix_nearly_id},  and 
 $\gamma_N$ is a constant 
 depending only on $d$ 
 and the matrix $\Sigma_N$.
To pass to the limit as $N\to+\infty$ in the ARE~\eqref{eq:ARE_nearly_id} we 
first note that~\eqref{eq:ARE_nearly_id} and~\eqref{eq:ARE_limit} are both AREs of the form~\eqref{eq:ARE} with
$$
{\cal R}_N={\cal R}=\,{\nu R\nu\over 2}\,
\qquad\qquad
{\cal Q}_N={A^T RA\over 2}\,+Q_N\,,
\qquad\qquad
{\cal Q}={A^T RA\over 2}\,+\hat Q\,,
$$
and the corresponding $2d\times 2d$ matrices of the form~\eqref{eq:inv_spaces} are given by
$$
{\cal H}_N=
\left(
\begin{array}{cc}
{\bf 0} & {\cal R}_N\\
{\cal Q}_N & {\bf 0}
\end{array}
\right)\,,
\qquad\qquad
{\cal H}=
\left(
\begin{array}{cc}
{\bf 0} & {\cal R}\\
{\cal Q} & {\bf 0}
\end{array}
\right)\,.
$$
Next, we claim 
 that the sequence $\Sigma_N$ is bounded w.r.t. the norm of the largest eigenvalue, which was introduced in~\eqref{eq:max_spec_norm}.
Indeed, 
by property {\it (iii)} of Proposition~\ref{prop:ARE}
\begin{equation}\label{speck}
\mathrm{spec}({\cal R}\Sigma_N)=\mathrm{spec}({\cal H}_N)\cap (0,+\infty)\,,
\end{equation}
and 
 the convergence ${\cal H}_N\to{\cal H}$ 
  implies the convergence of the 
   eigenvalues. Hence, for $N$ large enough
$$
\max\, \big\{\mathrm{spec}({\cal R}\Sigma_N)\big\} \leq \max\Big\{\mathrm{spec}({\cal H})\cap (0,+\infty)\Big\}+1
$$
and in particular the maximum eigenvalue of ${\cal R}\Sigma_N$ is bounded. Since ${\cal R}$ is symmetric positive definite, this implies that $\left\|\Sigma_N\right\|$ is bounded as well. Indeed, by denoting with $\lambda_{min}>0$ the smallest eigenvalue of ${\cal R}$, we have ${\cal R}-\lambda_{min}\Id_d\geq 0$ and 
$$
\left\|\Sigma_N\right\|=\,{1\over \lambda_{min}}\,\left\|\lambda_{min}\Sigma_N\right\|\leq
\max\mathrm{spec}({\cal R}\Sigma_N)<+\infty\,,
$$
thanks to Proposition~\ref{prop:matrices}--{\it (iii)}. 
Therefore, $\Sigma_N$ has a converging subsequence, that we denote with $\Sigma_{N_k}$, whose limit $\overline{\Sigma}$ satisfies~\eqref{eq:ARE_limit}, i.e., 
$$
\overline{\Sigma}\,{\nu R\nu\over 2}\,\overline{\Sigma}-\,{A^T RA\over 2}\,=\hat Q\,.
$$
If we can 
prove that $\overline{\Sigma}\in\Symp$, then we 
 have $\overline{\Sigma}=\Sigma$ by uniqueness in $\Symp$ of solutions to~\eqref{eq:ARE_limit}. 
 Since in general the limit of a sequence in $\Symp$ is only 
 semi--definite, 
 we further exploit the continuous dependence of the eigenvalues on 
 the coefficients of the matrix. 
 For $k$ large enough, by \eqref{speck},
$$
\min\, \big\{\mathrm{spec}({\cal R}\Sigma_{N_k})\big\} > \,{1\over 2}\,\min\Big\{\mathrm{spec}({\cal H})\cap (0,+\infty)\Big\}>0\,,
$$
and this implies
$$
\min\, \big\{\mathrm{spec}({\cal R}\overline{\Sigma})\big\} \geq \,{1\over 2}\,\min\Big\{\mathrm{spec}({\cal H})\cap (0,+\infty)\Big\}>0\,.
$$
Recalling again that ${\cal R}$ is symmetric positive definite, $\overline{\Sigma}$ cannot have zero as eigenvalue and this proves that the limit of $\Sigma_{N_k}$ is $\Sigma$.
Since we can repeat this argument to show that every subsequence of $\Sigma_N$ has a convergent subsequence whose limit is $\Sigma$, we conclude that $\Sigma_N\to\Sigma$ as $N\to+\infty$. 


Concerning the convergence of $\mu_N$ to $\mu$, we know that these vectors are respectively solutions to the linear systems
$$
{\cal B}'_N\mu_N=P'_N\,,
\qquad\qquad
{\cal B}_\infty\mu=P_\infty\,,
$$
where ${\cal B}'_N$ was defined in~\eqref{eq:weaker_invert},  ${\cal B}_\infty$ was defined in~\eqref{eq:limit_invert}, and the vector $P'_N,P_\infty$ are given, as in the previous sections, by
$$
P'_N=-Q_N H+(1-N)~{B_N\over 2}\,\Delta\,,
\qquad\qquad
P_\infty=-\hat Q H+\,{\hat B\over 2}\,\Delta\,.
$$
Here we use that 
 the matrix ${\cal B}_\infty$ is invertible by {\bf ($\mbox{U}_\infty$)}. Owing to ${\cal B}'_N\to{\cal B}_\infty$, it 
 follows that ${\cal B}'_N$ is invertible as well for $N$ large enough. In particular, for such $N$, 
$
\mu_N=({\cal B}'_N)^{-1} P'_N
$,
and we can pass to the limit as $N\to\infty$ proving that $\mu_N\to\mu$.


Finally, by passing to the limit in the explicit formulas
~\eqref{eq:KFP_matrix_nearly_id} and~\eqref{eq:lambda_nearly_id}
$$
\Lambda_N=R\big(\nu\Sigma_N+A\big)\,,\qquad\qquad \rho_N=-R\nu\Sigma_N\mu_N\,,
$$
$$
\lambda_N^i=\widetilde F^i_0-(\mu_N)^T\,{\Sigma_N \nu R \nu \Sigma_N\over 2}\, \mu_N+\Tr(\nu R\nu\Sigma_N+\nu RA)\,,
$$
it is easy to verify that $\Lambda_N\to\Lambda$, $\rho_N\to\rho$ and $\lambda_N^i\to\lambda$ for each $i$. 
Now we can pass 
 to the limit in the 
  formulas 
of the quadratic-Gaussian solutions~\eqref{formulaN} 
and deduce the convergence of the value function and of the invariant measure from the convergence of the coefficients, and this completes the proof.~~$\diamond$

\smallskip

\n {\it Proof of Theorem~\ref{thm:uniq}.} 
For any ${\frak m},{\frak n}\in{\cal P}_2(\R^d)$ 
\begin{align*}
\left(\hat V[{\frak m}]-\hat V[{\frak n}]\right)(x) &= \int_{\R^d} \left[(x-H)^T\,{\hat B\over 2}\,(\xi-\Delta)+(\xi-\Delta)^T\,{\hat B\over 2}\,(x-H)\right]\,d({\frak m}-{\frak n})(\xi)\\
&~~~~~+\int_{\R^d} (\xi-\Delta)^T \hat C (\xi-\Delta)\,d({\frak m}-{\frak n})(\xi)\\
&~~~~~+\left(\int_{\R^d} (\xi-\Delta)\,d{\frak m}(\xi)\right)^T\hat D\left(\int_{\R^d} (\xi-\Delta)\,d{\frak m}(\xi)\right)\\
&~~~~~-\left(\int_{\R^d} (\xi-\Delta)\,d{\frak n}(\xi)\right)^T\hat D\left(\int_{\R^d} (\xi-\Delta)\,d{\frak n}(\xi)\right)\,.
\end{align*}
Observe  that only the first term of this expression depends on the variable $x$. Then 
\begin{align*}
\int_{\R^d} \left(\hat V[{\frak m}]-\hat V[{\frak n}]\right)(x)\,& d\big({\frak m}-{\frak n}\big)(x)=\\
&~~~~~\int_{\R^d}\int_{\R^d} \left[(x-H)^T\,{\hat B\over 2}\,(\xi-\Delta)\right]\,d\big({\frak m}-{\frak n}\big)(x)\,d\big({\frak m}-{\frak n}\big)(\xi)\\
&~~~~~+\int_{\R^d}\int_{\R^d} \left[(\xi-\Delta)^T\,{\hat B\over 2}\,(x-H)\right]\,d\big({\frak m}-{\frak n}\big)(x)\,d\big({\frak m}-{\frak n}\big)(\xi)\,.
\end{align*}
The first term can be written as 
\begin{align*}
\int_{\R^d}\int_{\R^d} &\left[(x-H)^T\,{\hat B\over 2}\,(\xi-\Delta)\right]\,d\big({\frak m}-{\frak n}\big)(x)\,d\big({\frak m}-{\frak n}\big)(\xi)\\
&=\left(\int_{\R^d} x\,d\big({\frak m}-{\frak n}\big)(x)\right)^T\,{\hat B\over 2}\, \left(\int_{\R^d}\xi \,d\big({\frak m}-{\frak n}\big)(\xi)\right)\,.
\end{align*}
The second term can be treated in the same way, leading to
\begin{align*}
\int_{\R^d} \left(\hat V[{\frak m}]-\hat V[{\frak n}]\right)(x)\, d\big({\frak m}-{\frak n}\big)(x)&=\left(\int_{\R^d} x\,d\big({\frak m}-{\frak n}\big)(x)\right)^T\,\hat B\, \left(\int_{\R^d}x \,d\big({\frak m}-{\frak n}\big)(x)\right)
\end{align*}
It is now clear that if $\hat B$ is positive definite, then $\int_{\R^d} \left(\hat V[{\frak m}]-\hat V[{\frak n}]\right)(x)\, d({\frak m}-{\frak n})(x)\geq 0$ for all probability measures ${\frak m},{\frak n}$. Viceversa, for every fixed vector $\eta\in\R^d$, 
we can consider the multivariate Gaussian measures
${\frak m}={\cal N}(\eta,\Id_d)$ and ${\frak n}={\cal N}(0,\Id_d)$ on $\R^d$,
and apply the inequality above to obtain
$$
0\leq \int_{\R^d} \left(\hat V[{\frak m}]-\hat V[{\frak n}]\right)(x)\, d({\frak m}-{\frak n})(x)=\eta^T\,\hat B\,\eta\,,
$$
which implies $\hat B\geq 0$. This completes the proof of the equivalence 
 between monotonicity of the operator $\hat V$ and the positive semi-definitess of $\hat B$. The uniqueness of solutions when $\hat B\geq 0$ can be proved by 
  the same arguments of~\cite{Bardi, LL1, LL3}.~~$\diamond$

\begin{rem} \rm{ 
We believe that, 
 in general, the 
  assumptions~\eqref{eq:scale} and 
   {\bf ($\mbox{E}_\infty$)} 
   do not ensure condition {\bf ($\mathbf{E'}$)} with $Q=Q^N$, for $N$ large enough, mostly because the equation \eqref{eq:sylvester_nearlyid} 
    is equivalent to the symmetry of $\Sigma_N\nu R- RA$, 
    and this property is easily lost 
    under small perturbations.
When condition {\bf ($\mathbf{E'}$)} is violated, it would be interesting to verify whether the strategy resulting from the solution to~\eqref{eq:hyp_limit} could provide an $\ve$--Nash equilibrium for the $N$--player game, by using the techniques from~\cite{HCM:03,HCM:06,NCMH}. However, this investigation is beyond the scope of this paper and is left to a future work.
}\end{rem}

\begin{rem}  \rm{
An interesting open issue is  the rate of convergence of the solutions if the order of convergence of the data in~\eqref{eq:scale} is given. In the general case the rate of convergence of $\Sigma_N\to\Sigma$ does not follow from our analysis.
However an estimate can be obtained 
when we have an explicit formula for the solution of the Riccati equation, as  in the examples presented in Section~\ref{sec:explicit}, 
where it is possible to prove that
$
\|\Sigma-\Sigma_N\|={\cal O}(\|Q-Q_N\|).
$ 
 See also Remark 14 of \cite{Bardi}  on the possibility of expanding explicit solutions in powers of $1/N$.
}\end{rem}

\section[Examples]
{Examples}\label{sec:explicit}

So far we have used the abstract condition {\bf (E)} with {\bf (U)} 
 to translate the existence and uniqueness of solutions having the form~\eqref{eq:quad_gauss} to the system of PDEs~\eqref{eq:hjkfp} into algebraic matrix equations. 
In this section we show that in some 
  cases such conditions can 
 be easily verified and the solution to the PDEs~\eqref{eq:hjkfp} can be computed explicitly. For simplicity we limit ourselves to nearly identical players and the Mean-Field game. Therefore we focus on the corresponding conditions {\bf (E$'$)} and {\bf (U$'$)}.

In the last part we discuss some consensus models for which infinitely many solutions can be exhibited.
\subsection{Symmetric system}\label{symmetric}
Consider an $N$--players game with dynamics~\eqref{eq:sde} and costs~\eqref{eq:cost} and assume that {\bf (H)} holds, that players are nearly identical and 
 for all $i\in\{1,\ldots,N\}$ 
\begin{description}
\item{\it (a)} the dynamics~\eqref{eq:sde} involve 
drift  matrices $A^i\equiv A\in\Sym$ and diffusion matrices $\sigma=s\Id_d$ with $s\in\R\setminus\{0\}$;
\item{\it (b)} the matrix $R$ in the control costs~\eqref{eq:cost} satisfies $R=r\Id_d$ with $r>0$.
\end{description} 
Then it is easy to verify the part of {\bf (E$'$)} concerning solutions of~\eqref{eq:riccati_nearlyid} and~\eqref{eq:sylvester_nearlyid}. Indeed, 
both matrices 
$$
\nu=\,{s^2\over 2}\,\Id_d=: \bar\nu\Id_d
$$
and $R$ commute with any other matrix. Then, Sylvester's equation~\eqref{eq:sylvester_nearlyid} can be rewritten as
$$
r\big(A-A^T\big)=r\,\bar\nu\, (Y-Y)=0\,,
$$
i.e. it reduces to a symmetry condition on $A$, which is ensured by {\it (a)}. 
Moreover, an explicit expression of the matrix $\Sigma$ can be calculated. Indeed, the matrix ${2\over r}\,Q+A^2$ is symmetric and positive definite and thus admits a square root $E\in\Symp$, 
i.e.,
$$
E^2:= {2\over r}\,Q+A^2\,.
$$
If we now consider the ARE~\eqref{eq:ARE_nearly_id}, we find that
$$
r\,{\bar\nu^2\over 2}\,\Sigma^2=\,Q+\, {r\over 2}\,A^2=\, {r\over 2}\, E^2\,,
$$
which implies
$$
\Sigma=\,{1\over \bar\nu}\, E =  {1\over \bar\nu}\sqrt{{2\over r}\, Q+A^2}\,.
$$

 To verify the part of condition {\bf (E$'$)} dealing with the matrix ${\cal B}'$, we assume in addition, using the notations
of Lemma~\ref{lem:symmetry}, 
 \begin{description}
\item{\it (c)} the primary costs of displacement $B_i$ in $Q^i$ satisfy $B_i=B\geq 0$  for all $i=1,\ldots,N$.
\end{description} 
Then we can rewrite~\eqref{eq:weaker_invert} as
$$
{\cal B}'\,= \,Q+\,\frac{r}{2}\,A^2+\frac{N-1}2B\, =  \,\frac{r}{2}\, E^2+\frac{N-1}2B\,,
$$
so ${\cal B}'$ is invertible and both {\bf (E$'$)} and  {\bf (U$'$)} are satisfied. 
Thus the linear system~\eqref{eq:mu_nearly_id} has the 
unique solution 
$$
\mu = 
({\cal B}')^{-1}\left(Q H+\frac{N-1}2B\Delta\right)\,,
$$
where $H$ and $\Delta$ are 
the reference positions of the players, as in Definition~\ref{def:nearly_id}. 
Now the expressions found for $\Sigma$ and $\mu$ can be used in \eqref{eq:KFP_matrix_nearly_id} and \eqref{eq:lambda_nearly_id} to obtain $\Lambda,\rho$ and $\lambda^1,\ldots\lambda^N$, completing the explicit construction of the unique solution of quadratic--Gaussian type.

 In conclusion, for games with $N$ nearly identical players which satisfy {\bf (H)}, conditions {\it (a)}, {\it (b)} and {\it (c)} are sufficient to guarantee the existence of a unique solution to~\eqref{eq:hjkfp} of the form~\eqref{eq:ansatz_nearly_id}, and hence of a unique affine Nash equilibrium strategy given by
$$
\overline{\alpha}(x) = Ax+E(x-\mu) 
\,.
$$

For the large population limit we assume that the scaled coefficients satisfy~\eqref{eq:scale} as $N\to+\infty$
and
\begin{description}
\item{\it ($\mbox{a}_\infty$)} $A\in\Sym$, $\nu= \bar\nu\Id_d$ with $\bar \nu>0$;
\item{\it ($\mbox{b}_\infty$)} $R=r\Id_d$ with $r>0$;
\item{\it ($\mbox{c}_\infty$)}  $\hat B\geq
 0$.
\end{description} 
Once again, Sylvester's equation reduces to the symmetry of $A$, the solution to ARE~\eqref{eq:ARE_limit} can be given explicitly as $\Sigma=\,{1\over \bar\nu}\, \sqrt{{2\over r}\,\hat Q+A^2}$, and the invertibility of ${\cal B}_\infty\,= \,\hat Q+\,\frac{r}{2}\,A^2+\frac{\hat B}{2}$ is immediate. Moreover, it is easy to verify that $\Sigma_N\to\Sigma$ and that $\mu_N\to\mu$, which in turn imply convergence of the unique solution to the $N$--players game to the unique solution of the Men-Field game~\eqref{eq:hjkfp_limit}.

\begin{rem} \rm{
The equivalence between $A$ symmetric and Sylvester's equation implies 
that 
if $A$ is not symmetric condition {\bf (E$'$)} 
 fails and therefore no solution of~\eqref{eq:hjkfp} with a quadratic value function exists. Indeed, 
 in this case one gets an affine vector $q(x)=\Lambda x+\rho$ which solves, together with a multivariate Gaussian $m$, the equation
$$
-\Tr(\nu\,\D^2 m)-\displaystyle\Div\Big(m {\partial H\over\partial p}(x,q)\Big)=0\,,
$$
but $q$ is 
not the gradient of a 
 quadratic 
  function of the form~\eqref{eq:ansatz_nearly_id}, because $\Lambda\notin\Sym$. 
}\end{rem}

\subsection{Non--defective system}\label{non-def}
In this section we extend 
the previous analysis 
beyond the symmetry assumption on the drift matrix $A$ to the case of 
$A$ 
 \emph{non--defective}. We recall that a matrix $M\in\Mat_{d\times d}(\R)$ is said to be non--defective if, for every eigenvalue $\lambda\in\mathrm{spec}(M)$, the corresponding eigenspace has dimension equal to the multiplicity of $\lambda$ or, equivalently, when there exists a base of $\R^d$ consisting of right (or left) eigenvectors of $M$. 
%
%
\begin{prop}\label{prop:symmetrizers} Let $M$ be any $d\times d$ real matrix. Then the following properties hold.
\begin{description}
\item{\it (i)} There exists an invertible and symmetric symmetrizer for $M$, i.e. there exists a matrix $Y\in\Mat_{d\times d}(\R)$ such that
$$
\det(Y)\neq 0\,,
\qquad\qquad
Y^T=Y\,,
\qquad\qquad
YM=M^TY\,.
$$
\item{\it (ii)} If $M$ is non--defective, then the symmetrizer $Y$ can be chosen positive definite.
\item{\it (iii)} If $M$ is non--defective, $\sigma$ is invertible and we consider a linear stochastic differential equation
\begin{equation}\label{eq:sde_nondef}
dx_t=(M x_t-\alpha_t)dt+\sigma dW_t\,,
\end{equation}
then there exists a linear change of coordinates $x\mapsto \xi$ such that~\eqref{eq:sde_nondef} can be rewritten in the form
\begin{equation}\label{eq:sde_symm}
d\xi_t=(\widetilde M\xi_t-\tilde\alpha_t)dt+\tilde\sigma dW_t\,,
\end{equation}
with $\widetilde M$ symmetric matrix and $\tilde\sigma$ invertible.
\end{description}
\end{prop}

\n{\it Sketch of the proof.} We refer to~\cite{TaZa} for the proof of {\it (i)}.
{\it (ii)} Bhaskar proved in~\cite{Bhas} that the real symmetrizer $Y$ can be chosen of the explicit form $Y=(U^{-1})^TV^T$, where $U$ and $V$ are $d\times d$ (complex) matrices having as columns left and right eigenvectors of $A$, respectively. By orthogonality between left and right eigenvectors, it is immediate to see that $V^TU=\Id$, and hence that
$$
Y=(U^{-1})^TV^T=V V^T
$$
is real and positive definite. \\
{\it (iii)} By choosing $P\in\mathrm{Mat}_{d\times d}(\R)$ orthogonal such that
$$
Y=P^T D P\,,
$$
with $D$ diagonal, and $Z\in\Symp$ such that $D=Z^2$, one can verify that the linear change of coordinates $\xi=ZPx=Z^{-1}PYx$ allows to rewrite the stochastic linear equation in the required form~\eqref{eq:sde_symm} with
$$
\widetilde M:= Z^{-1}PYMP^TZ^{-1}\,,
\qquad\qquad 
\tilde \alpha := ZP\alpha\,,
\qquad\qquad 
\tilde\sigma:= ZP\sigma\,.
$$
This completes the proof.~~$\diamond$


Observe that, if we consider a differential game with $N$ nearly identical players, dynamics~\eqref{eq:sde} and costs~\eqref{eq:cost} and if we assume that the drift matrix $A$ is non--defective, then we can perform the change of coordinate in Proposition~\ref{prop:symmetrizers} to pass to a new SDE with symmetric drift. Let us denote, as in the proof of the proposition, with $Y=P^TZ^2P$ the symmetrizer matrix for $M$ with $P$ orthogonal matrix and $Z$ diagonal and positive definite matrix. Then, the new game will have costs given by
$$
\widetilde J^i(\Xi,\tilde\alpha^1,\ldots,\tilde\alpha^N):= \liminf_{T\to\infty}\,{1\over T}\,\vatt\left[
\int_0^T \,{(\tilde\alpha_t^i)^T\widetilde R\tilde\alpha_t^i\over 2}\,+\sum_{j,k=1}^N (\Xi^j-\overline{\Xi_i}^j)^T \widetilde Q^i_{jk}(\Xi^k-\overline{\Xi_i}^k)\,dt
\right]\,,
$$
where the new variables $\Xi\in\R^{Nd}$ and $\tilde\alpha^1,\ldots,\tilde\alpha^N\in\R^d$ satisfy for $k\in\{1,\ldots,N\}$
$$
\Xi^k=ZPX^k\,,\qquad \qquad\qquad \tilde\alpha^k=ZP\alpha^k\,,
$$
and the matrices $\widetilde R$ and $\widetilde Q^i_{jk}$ are given by
$$
\widetilde R=Z^{-1}PRP^{-1}Z^{-1}\,,\qquad \qquad\qquad \widetilde Q^i_{jk}=Z^{-1}PQ^i_{jk}P^{-1}Z^{-1}\,,
$$
for the same matrices $Z\in\Symp$ and $P$ orthogonal used to symmetrize $A$. One can easily verify that replacing {\it (a)} and {\it (b)} with
\begin{description}
\item{\it (a$'$)} dynamics~\eqref{eq:sde} are given by drift matrices $A^i\equiv A$ with $A$ non--defective and diffusion matrices $\sigma=s P^TZ^{-1}$ with $s\in\R\setminus\{0\}$;
\item{\it (b$'$)} matrix $R$ in control costs~\eqref{eq:cost} satisfies $R=rY$ with $r>0$;
\end{description} 
it is possible to repeat the 
 arguments of Sect.~\ref{symmetric}, after the change of coordinates $\xi=ZPx$, and to prove that for games with $N$ nearly identical players 
 satisfying {\bf (H)}, 
  {\it (a$'$)}, {\it (b$'$)} and {\it (c)} 
  there exists 
   a unique solution to~\eqref{eq:hjkfp} of the form~\eqref{eq:ansatz_nearly_id}. Similar conditions in the case of $A$ only non--defective can be given for the limit problem~\eqref{eq:hjkfp_limit} as well.

\begin{rem} \rm{
These arguments apply to general $N$--players games too, without the assumption of nearly identical players, whenever all drift matrices $A^1,\ldots, A^N$ are symmetric or simultaneously symmetrizable (i.e., if the symmetrizers $Y^1,\ldots,Y^N$ given by Proposition~\ref{prop:symmetrizers} coincide). For more general games one has either to require that more blocks of the cost matrices $Q^i$ are null or to study the $Nd$ dimensional linear SDE for $X_t=(X^1_t,\ldots,X^N_t)$.
}\end{rem}

\subsection{Consensus models}\label{sec:consensus}
In this section we apply the previous theory to some simple models of consensus in multi-agent systems inspired by \cite{NCMH} to which we refer for motivations and bibliography. Consider costs whose part depending on the state is 
\begin{equation}\label{eq:special_cost}
F^i(X^1,\ldots,X^N)=\frac 1{N-1}\sum_{j\neq i} (X^i-X^j)^T P^N(X^i-X^j) , \qquad i=1,\dots, N ,
\end{equation}
with $P^N\in\Symp$.  
Then each player seeks a position as close as possible to the positions of the other players. Note that $F^i$ has the symmetry property {\bf (S)} and the blocks of the matrix $Q^i$ are the following, in the notations of Section \ref{sec:limit},
\begin{equation}\label{eq:special_blocks}
Q^N=P^N,\quad B^N=-\frac 2{N-1}P^N , \quad C_i^N = \frac 1{N-1}P^N ,\quad D_i^N=0 ,
\end{equation}
and reference states $H=\Delta=0$.

Assume also that the dynamics and cost of the control is the same for all players, so that they are nearly identical.
For simplicity suppose also the conditions \emph {(a)} and \emph {(b)} of Section~\ref{symmetric}, although the analysis can be carried over to merely non-defective matrices $A$. 

The matrix ${\cal B}'$ defined in \eqref{eq:weaker_invert} is
$$
{\cal B}'= P^N+{r\over 2}A^2\,+(N-1)~{B^N\over 2}={r\over 2}A^2 \qquad \forall\,N ,
$$
and the equation~\eqref{eq:mu_nearly_id} for the mean $\mu\in\R^d$ of the distribution of the players is ${\cal B}'\mu=0$. Then the condition {\bf (E$'$)}  holds and the existence part of Theorem \ref{thm:NIDplay} 
proves the following:

\noindent \emph{there exists an identically distributed quadratic-Gaussian solution $(v_N,m_N,\lambda_N)$ with  $m_N \sim{\cal N}(\mu,\Sigma_N^{-1})$ 
if and only if  $\mu$ is such that $A\mu=0$}. 

Moreover, the covariance matrix is the same for all $\mu$
\begin{equation}\label{eq:special_sigma}
\Sigma_N = {1\over \bar\nu}\sqrt{{2\over r}\, P^N+A^2}\,,
\end{equation}
and the Nash feedback equilibrium is
\begin{equation*}\label{eq:special_alfa}
\overline\alpha_N(x) = Ax + {1\over \bar\nu}\sqrt{{2\over r}\, P^N+A^2}\,(x-\mu)\,.
\end{equation*}
In particular, there is a unique solution if and only if $\det A\ne 0$ and then the mean 
 is $\mu=0$. 
 
 The large population limit is straightforward if we assume that $P^N\to \hat P>0$ as $N\to\infty$. Condition \eqref{eq:scale} is satisfied with 
 \[
  \hat Q = \hat P
\qquad
 \hat B = -2 \hat P
\qquad
 \hat C = \hat P
\qquad
 \hat D =0
 \]
 and passing to the limit in \eqref{eq:special_sigma} one gets $\hat\Sigma = {1\over \bar\nu}\sqrt{{2\over r}\, \hat P+A^2}$. Then \emph{the Mean-Field Games PDEs \eqref{eq:hjkfp_limit} have  a solution $(v,m,\lambda)$ with $m \sim{\cal N}(\mu,\hat\Sigma^{-1})$ if and only if $\mu\in\R^d$ is such that $A\mu=0$. }
 
  \begin{rem} \rm{
 The paper by Nourian et al. \cite{NCMH} uses  a cost term of the form 
 \[
 F^i(X^1,\ldots,X^N):=\left
 | X_i - \frac 1{N-1}\sum_{j\neq i} X^j\right
 |^2 , \qquad i=1,\dots, N ,
 \]
 instead of \eqref{eq:special_cost}. Then 
  $Q^N=\Id_d$ and  $B^N=-\frac 2{N-1}\Id_d$. Although the secondary costs are different, this is a special case of the above, with $P^N=\Id_d$ for all $N$. In fact ${\cal B}'$ is the same and $\Sigma_N=\hat\Sigma = {1\over \bar\nu}\sqrt{{2\over r}\, \Id_d+A^2}$, 
  so we get the same conclusions on the quadratic-Gaussian solutions of $N$-person as well as Mean-Field game.
  }
 \end{rem}
 \begin{rem} \rm{The existence of infinitely many Gaussian solutions in a Mean-Field Game model of population distribution with rewarding imitation among players 
  was first observed by Gu\'eant \cite{G:09}. 
 In \cite{Bardi}   the   
  LQ Mean-Field Game of Section \ref{sec:limit} was studied for $d=1$,  $A=0$, and $H=\Delta$ and it was observed that for $\hat Q\ne -\hat B/2$ there is a unique quadratic-Gaussian solution with $\mu =H$, whereas for $\hat Q = -\hat B/2$ there are infinitely many, one for any $\mu\in\R^d$.
   }
 \end{rem}

\begin{rem} \rm{
Analogous computations can be made also for consensus models where the dynamics and cost of the control are not the same for all players. In such a case, the matrix ${\cal B}'$ above is replaced by the matrix ${\cal B}$ defined in~\eqref{eq:matrixB} which becomes here
\begin{equation}\label{eq:rank1_pert}
{\cal B}=\,{1\over N-1}\left(
\begin{array}{c}
\sqrt{P}\\
\vdots\\
\sqrt{P}
\end{array}
\right)\left(
\begin{array}{c}
\sqrt{P}\\
\vdots\\
\sqrt{P}
\end{array}
\right)^T\!\!\!-\,\mathrm{diag}\! \left({N\, P\over N-1}\,+\,{(A^1)^T R^1 A^1\over 2}\,,\ldots,{N\, P\over N-1}\,+\,{(A^N)^T R^N A^N\over 2}\,\right)
\end{equation}
i.e., ${\cal B}$ a ``block--rank--one'' perturbation of a $Nd\times Nd$ block--diagonal matrix, with $d\times d$ blocks $-{N\, P\over N-1}\,-\,{(A^\alpha)^T R^\alpha A^\alpha\over 2}$. It is easy to verify that, if there exists $\xi\neq 0$ in $\bigcap_{\alpha=1}^N\mathrm{Ker}(A^\alpha)$, then $\eta:=(\xi,\ldots,\xi)^T\in\R^{Nd}$ provides a solution to ${\cal B}\eta = 0$, which is~\eqref{eq:system_mi} for this model. 
Therefore, in this case we find infinitely many Gaussians solutions $m^i\sim{\cal N}(\mu^i,(\Sigma^i)^{-1})$ for each player $i=1,\ldots,N$, with $\mu^i=\xi$ and $\Sigma^i$ solving~\eqref{eq:cond1}. 
In fact, owing to known formulas for matrices of the form~\eqref{eq:rank1_pert}, it is possible to prove that the
invertibility of ${\cal B}$ is equivalent to $\bigcap_{\alpha=1}^N\mathrm{Ker}(A^\alpha)=\{0\}$, providing a complete understanding of the conditions required to have a unique QG solution for this consensus model.
} \end{rem}

\smallskip

\n{\bf Acknowledgements.} We thank the anonymous referees for several useful suggestions. Part of this research was done while the second-named author was holding post-doc positions at the Departments of Mathematics of the Universities of Padova and Roma-Tor Vergata.

\end{document}